\documentclass[10pt]{article}
\usepackage{mathrsfs}
\usepackage[centertags]{amsmath}
\usepackage{amsfonts}
\usepackage{amssymb}
\usepackage{amsthm}
\usepackage{epsfig}
\usepackage{color}
\allowdisplaybreaks[4]
\usepackage{indentfirst}
\usepackage{amssymb,color}
\allowdisplaybreaks
\definecolor{c20}{rgb}{1.00,0.00,0.00}
\definecolor{c30}{rgb}{0.,0.,1.}
\definecolor{c40}{rgb}{1,0.1,0.7}
\definecolor{c50}{rgb}{1,0,0}
\definecolor{c60}{rgb}{1,0.9,0.1}
\definecolor{c70}{rgb}{0.50,1.00,0.00}
\definecolor{c80}{rgb}{0.00,1.00,0.00}

\date{}

\newtheorem{theorem}{Theorem}[]
\newtheorem{proposition}{Proposition}[]

\newtheorem{lemma}{Lemma}[]
\newtheorem{remark}{Remark}[]
\newtheorem{definition}{Definition}[]
\newtheorem{example}{Example}[]
\numberwithin{equation}{section}

\def\R{\operatorname*{\mathbb{R}}}

\def\CTE{\operatorname*{\mathrm{CTE}}}
\def\VaR{\operatorname*{\mathrm{VaR}}}

\def\C{\operatorname*{\mathcal{C}}}
\def\L{\operatorname*{\mathcal{L}}}
\def\S{\operatorname*{\mathcal{S}}}
\begin{document}
\title{Second-order asymptotics for convolution of distributions \\ with light tails}
\author{{  Zuoxiang Peng\quad Xin Liao\thanks{Corresponding author. Email: liaoxin2010@163.com}}\\
{\small\it School of Mathematics and Statistics, Southwest
University, 400715 Chongqing, China }}
\maketitle

\begin{quote}
{\bf Abstract}~~In this paper, asymptotic behaviors of convolutions of
distributions belonging to two classes of distributions with
light tails are considered, respectively. The precise
second-order tail asymptotics of the convolutions are derived under
the condition of second-order regular variation.



\end{quote}

\section{Introduction}
\label{sec1}

Let $X$ and $Y$ be two independent and nonnegative random variables with
distribution $F$ and $G$, respectively. The distribution of the sum $X+Y$,
written as $F*G$, is called the convolution of $F$ and $G$.

A distribution $F$ on $[0,\infty)$ is said to belong to the class $\L_{\alpha}$ for some $\alpha\geq 0$,
if its right tail satisfies $\overline{F}(t)=1-F(t)>0$ for all $t> 0$ and
the relation
\begin{equation}\label{eq1.1}
\lim_{t\to \infty}\frac{\overline{F}(t-u)}{\overline{F}(t)}=e^{\alpha u}
\end{equation}
holds for all $u$. For $\alpha=0$, the
class $\L_{0}$ reduces to the well-known class of long-tailed
distributions. Clearly, the class $\L_{\alpha}$ is related to the
class $RV_{-\alpha}$ of regularly varying functions with exponent
$-\alpha$ by the fact that
\begin{eqnarray*}
F\in {\L}_{\alpha} \quad \text{if and only if } \quad \overline{F}(\ln t)\in RV_{-\alpha}.
\end{eqnarray*}
Applying Karamata's representation theorem for regularly varying
functions (see Bingham et al. (1987) and Cline (1986)), we know that $F\in \L_{\alpha}$ if and only
if $\overline{F}(t)$ can be expressed as
\begin{equation}\label{eq1.2}
\overline{F}(t)=c(t)\exp\left( -\int_{0}^{t} \alpha(y)dy \right)
\end{equation}
with $\lim_{t \to \infty}c(t)=c>0$ and $\lim_{t\to \infty}\alpha(t)
=\alpha$.

An important subclass of $\L_{\alpha}$ is the class of convolution
equivalent distributions $\S_{\alpha}$. We say that $F$ belongs to
the class $S_{\alpha}$ if $F\in \L_{\alpha}$ and the limit
\begin{equation}\label{eq1.3}
\lim_{t\to \infty} \frac{\overline{F*F}(t)}{\overline{F}(t)}=2m_{F}(\alpha)
\end{equation}
exist and is finite, where the constant
$m_{F}(\alpha)=\int_{0}^{\infty} e^{\alpha u}dF(u)$; see Chistyakov
(1964) and Chover et al. (1973a, 1973b), Cline (1987). The class $\S:= \S_{0}$ is called the class of
subexponential distributions.

Properties of those mentioned classes have been extensively
investigated by many researchers and have been applied to many
fields. Embrechts and Goldie (1980) proved that $\L_{\alpha}$ is
closed under convolution, and gave sufficient conditions for the
statements $F\in \S \Leftrightarrow F*G\in \S$ and $F*G\in \S
\Leftrightarrow pF+(1-p)G\in \S$ with some (all) $p\in (0,1)$. Some
more results on $\S_\alpha$ are presented by Embrechts and Goldie
(1982). For $F, G\in \L_{\alpha}$, Cline (1986) investigated the
relationship between $\overline{F*G}$ and its components
$\overline{F}$ and $\overline{G}$. Further, Cline (1987) derived
some results on distribution tails of random sums from
$\S_{\alpha}$, and gave more closure and factorization properties
for distributions belonging to $\L_{\alpha}$ and $\S_{\alpha}$.
Recently, Zachary and Foss (2006) derived the asymptotic of tail distribution of convolution of distributions belonging to $\S_{\alpha}$, and further gave an distributional asymptotic for the supremum of a random walk with increments in $\S_{\alpha}$.
By using $h$-insensitivity function, more properties of
convolutions of long-tailed and subexponential distributions were investigated by Foss et al. (2009).
For more studies on random sums, random walk and queue theory related to the classes $\L_{\alpha}$ and
$\S_{\alpha}$, we refer to Albin (2008), Borovkov (1976), Foss (2007), Pakes (2007), Shimura and Watanabe (2005), Zachary (2004).

With motivation from Cline (1986), in this paper we are interested
in the second-order tail asymptotic expansions of convolutions of
the following two classes of distributions:
\begin{equation}\label{eq1.4}
\overline{F}(t)=e^{-\alpha t+ \chi(t)}, \quad \chi(t)\in
RV_{\rho},\, 0<\rho<1
\end{equation}
and
\begin{equation}\label{eq1.5}
\overline{F}(t)=b(t)e^{-\alpha t}, \overline{G}(t)=c(t) e^{-\alpha t}, \quad
b(t)\in RV_{\beta}, c(t)\in RV_{\gamma},\, \beta, \gamma \in \R.
\end{equation}
Clearly, $F$ and $G$ given by \eqref{eq1.4} and \eqref{eq1.5} not only belong to $\L_{\alpha}$, but also are Weibull-type distributions.
For such distributions $F$ and $G$, the first-order tail asymptotics of $F*G$ have been derived by Cline (1986).

In recent literature, more and more researchers focus on the
second-order asymptotic behaviors for the sake of understanding
precisely the tail behaviors of risks, ruin probabilities and random
summation. Hua and Joe (2011) obtained the second-order
approximation of the conditional tail expection $\CTE_{p}(X)$ for
risk $X$ with its survival function having the property of the
second-order regular variation. Degen et al. (2010) and Mao and Hu (2013)
derived the second-order approximations
of the risk concentrations $C_{\VaR}(p)$ and $C_{\CTE}(p)$,
respectively. Baltr\={u}nas (2005) investigated second-order
behaviour of ruin probabilities with subexponential claim-size. Lin
(2012) established second-order asymptotics for ruin probabilities
in a renewal risk model with heavy-tailed claims. For the
convolutions of the regularly varying distributions, the
second-order or higher-order tail asymptotics have been discussed by
Borovkov and Borovkov (2001),  Geluk et
al. (1997), Hashorva et al. (2014) and among others. Barbe and
McCormick (2008) and Lin (2014) respectively established the second-order tail
asymptotics of convolution for rapidly varying subexponential
distributions and second order subexponential distributions.
For related works on the higher-order asymptotics of random
sum, see, e.g., Geluk (1996),
Omey and Willekens (1986, 1987), Albrecher et al. (2010). To the best of our knowledge, there are no studies
on higher-order expansions for convolution of
distributions belonging to $\L_{\alpha}$ with $\alpha>0$. The main goal of this paper is to investigate the
second-order asymptotics of convolution of light tail
distributions $F, G$ which satisfied \eqref{eq1.4} and
\eqref{eq1.5}. In order to get the desired results, we assume that
$\chi(t)$, $b(t)$, $c(t)$ in \eqref{eq1.4} and \eqref{eq1.5} are
second-order regularly varying functions.

The rest of this paper is organized as follows. In Section \ref{sec2},
some preliminary concepts and results of second-order regularly varying functions are presented, which will be
used to prove the main results.
The main results and some illustrating examples are given in Section \ref{sec3}. All proofs are deferred to Section \ref{sec4}.

\section{Preliminaries}
\label{sec2}

For analysis on tail asymptotics of convolutions of $F$ and $G$ satisfied \eqref{eq1.4} and \eqref{eq1.5},
the theory of regularly variation on survivor functions plays an important role.
We refer to de Haan and Ferreira (2006) and Resnick (1987) for standard references on regular variation.

\begin{definition}
A measurable function $\chi: {\R}_{+}\to \R$ that is eventually positive is
regularly varying at $\infty$ with some $\alpha \in \R$
(written $\chi\in RV_{\alpha}$) if for any $u>0$,
\begin{equation}\label{eq5.1}
\lim_{t\to \infty} \frac{\chi(tu)}{\chi(t)}=u^{\alpha}.
\end{equation}
We call $\alpha$ the index of variation. If $\alpha=0$, $\chi(t)$ is
said to be slowly varying at $\infty$.
\end{definition}

The following result is the famous Karamata's Theorem.

\begin{lemma}\label{lem2.1}
\emph{(de Haan and Ferreira, 2006, Theorem B.1.5; Resnick, 1987, Karamata's Theorem 0.6)}
Suppose that $\chi\in RV_{\alpha}$ for some $\alpha\in \R$ and $t_{0}>0$ such that $\chi(t)$ is positive for $t>t_{0}$.
If $\alpha\geq -1$, then
\begin{equation}\label{eq5.2}
\lim_{t\to \infty} \frac{t\chi(t)}{\int_{t_{0}}^{t}\chi(s)ds}=\alpha+1.
\end{equation}
If $\alpha<-1$, or $\alpha=-1$ and $\int_{0}^{\infty} \chi(s)ds<\infty$, then
\begin{equation}\label{eq5.3}
\lim_{t\to \infty} \frac{t\chi(t)}{\int_{t}^{\infty}\chi(s)ds}=-\alpha-1.
\end{equation}
Conversely, if \eqref{eq5.2} holds with $-1<\alpha<\infty$, then $\chi\in RV_{\alpha}$;
if \eqref{eq5.3} holds with $-\infty<\alpha<-1$, then $\chi\in RV_{\alpha}$.
\end{lemma}

Karamata's Theorem described the effect of integration on a regularly varying function.
When an $\alpha$-varying function is differentiated, the associated property was investigated by Proposition 0.7 of Resnick (1987) which is cited as follows.

\begin{lemma}\label{lem2.2}
\emph{(Resnick, 1987, Proposition 0.7)}
Suppose $\chi: {\R}_{+}\to {\R}_{+}$ is absolutely continuous with density $\chi'$ so that
$\chi(x)=\int_{0}^{x} \chi'(t)dt.$
If
\begin{equation}\label{eq5.4}
\lim_{t\to \infty} \frac{t\chi'(t)}{\chi(t)}=\alpha,
\end{equation}
then $\chi\in RV_{\alpha}$. Conversely, if $\chi\in RV_{\alpha}$,
$\alpha\in \R$, and $\chi'$ is monotone then \eqref{eq5.4} holds and if $\alpha\neq 0$, then $(\emph{sgn} \alpha)\chi'(x)\in RV_{\alpha-1}$.
\end{lemma}

The following definition of the second-order regular variation comes
from de Haan and Ferreira (2006) and Geluk et al. (1997).

\begin{definition}\label{eq5.5}
A measurable function $\chi: {\R}_{+}\to \R$ that is eventually positive is
said to be of second-order regularly variation with the first-order parameter $\alpha \in \R$
and the second-order parameter $\rho\leq 0$, denoted by $\chi \in 2RV_{\alpha, \rho}$,
if there exists some ultimately positive or negative function $A(t)$ with $\lim_{t\to \infty}$ A(t)=0 such that
\begin{equation}
\lim_{t\to \infty} \frac{\frac{\chi(tx)}{\chi(t)}-x^{\alpha}}{A(t)}
=x^{\alpha}\frac{x^{\rho}-1}{\rho}, \quad x>0.
\end{equation}
Here, $\frac{x^{\rho}-1}{\rho}$ is interpreted as $\log x$ when $\rho=0$, $A(t)$ is referred as auxiliary function of $\chi$,
and $\rho$ governs the speed of convergence in \eqref{eq5.1}.
\end{definition}

Next result concerns the Drees-type inequalities for $RV$ functions and $2RV$ functions
which establishes uniform inequalities.

\begin{lemma}\label{lem2.3}
\emph{(de Haan and Ferreira, 2006, Theorem B.1.10, Theorem 2.3.9; Drees, 1998)}
If $\chi\in RV_{\alpha}$ with $\alpha \in \R$, for each $\epsilon, \delta >0$, there is a $t_{0}=t_{0}(\epsilon, \delta)$ such that for $t, tx \geq 0$,
\begin{equation}\label{eq3.0}
\left| \frac{\chi(tx)}{\chi(t)}-x^{\alpha} \right| \leq \varepsilon x^{\alpha}\max\left( x^{\delta}, x^{-\delta} \right).
\end{equation}
Further, if $\chi \in 2RV_{\alpha, \rho}$ with auxiliary function $A(t)$ and $\rho\leq 0$,
then for any $\varepsilon, \delta >0$, there exist an auxiliary function $A_{1}(t)$, $A_{1}(t)\sim A(t)$ as $t \to \infty$, and $t_{0}=t_{0}(\varepsilon, \delta)>0$ such that for all $t, tx>t_{0}$,
\begin{equation}\label{eq3.1}
\left| \frac{\frac{\chi(tx)}{\chi(t)}-x^{\alpha}}{A_{1}(t)}
-x^{\alpha}\frac{x^{\rho}-1}{\rho} \right|
\leq \varepsilon x^{\alpha+\rho}\max\left( x^{\delta},x^{-\delta} \right).
\end{equation}
\end{lemma}

To end this section, we provide the following nice representation of
$2RV_{\alpha,\rho}$ with $\rho<0$ given by Hua and Joe (2011).

\begin{lemma}\label{lem2.4}
\emph{(Hua and Joe, 2011)}
Let $\alpha \in \R$, $\rho<0$ and $|A(t)|\in RV_{\rho}$. Then $\chi\in 2RV_{\alpha, \rho}$
with auxiliary function $A(t)$ if and only if there exists a constant $a>0$ such that
\[\chi(t)=at^{\alpha}\left( 1+\frac{1}{\rho}A(t) +o\left( A(t) \right) \right) \]
as $t\to \infty$.
\end{lemma}

\section{Main results}
\label{sec3}

In this section, we provide the main results. For
$\overline{F}(t)=e^{-\alpha t+\chi(t)}$ given by \eqref{eq1.4}, the
second-order tail asymptotic of $\overline{F*F}$ are presented in
the following theorem by assuming that $\chi(t)\in 2RV_{\rho,
\rho_{1}}$, $0<\rho<1$, $\rho+ \rho_{1}< 0$.

\begin{theorem}\label{th1}
Let $\overline{F}(t)=e^{-\alpha t+\chi(t)}$ with $\alpha>0$ and
$\chi(t)$ is eventually differentiable such that $\chi'(t)$ is
nonincreasing eventually. Assume that $\chi(t)\in
2RV_{\rho,\rho_{1}}$ with $0<\rho<1$, $\rho+\rho_{1}<0$ and
auxiliary function $A_{1}(t)$. Then, for large $t$ we have
\begin{equation}\label{eq2.1}
\frac{\chi^{\frac{1}{2}}\left( \frac{t}{2} \right)\overline{F*F}(t)}
{t\overline{F}^{2}\left( \frac{t}{2} \right)}
= \frac{\alpha}{2}\sqrt{\frac{\pi}{\rho(1-\rho)}}
-\frac{2^{\rho-5}\alpha\sqrt{\pi}(\rho-2)(\rho-3)}{(\rho(1-\rho))^{\frac{3}{2}}\chi(t)}
-\frac{1}{2^{\rho}}\sqrt{\frac{\rho\pi}{1-\rho}}\frac{\chi(t)}{t}
+ o\left( A_{1}(t)\chi(t)+\frac{1}{\chi(t)}+ \frac{\chi^{\frac{1}{2}}(t)}{t} \right).
\end{equation}
\end{theorem}

\begin{remark}
For $\overline{F}(t)=e^{-\alpha t+\chi(t)}$, we have $m_{F}(\alpha)=\infty$ which
implies that $F\in \L_{\alpha}$ but $F \not \in \S_{\alpha}$.
\end{remark}

\begin{example}
Let risk $X$  have distribution tail $\overline{F}(x)=e^{-\alpha x
+\chi(x)}$ with $\chi(t)=\alpha b t^{\frac{1}{2}} \left(
1+\frac{1}{2}A_{1}(t) +o(A_{1}(t)) \right)$, $b>0$, $|A_{1}(t)|\in
RV_{-2}$. Then by Lemma \ref{lem2.4}, $\chi(t)\in
2RV_{\frac{1}{2},-2}$, and the second-order expansion \eqref{eq2.1} leads to
\[\overline{F*F}(t)
=t^{\frac{3}{4}}\overline{F}^{2}\left( \frac{t}{2} \right)
\left( 2^{\frac{1}{4}}\left( \frac{\alpha \pi}{b} \right)^{\frac{1}{2}}
- 2^{-\frac{1}{4}}\left(  \frac{\pi}{\alpha b} \right)^{\frac{1}{2}}\left( 1+\frac{15}{8\alpha b^{2}} \right)t^{-\frac{1}{2}}+o\left( t^{-\frac{1}{2}} \right) \right)\]
for large $t$.
\end{example}

Theorem \ref{th4} provides the second-order expansions of
$\overline{F*G}$ where $F$ and $G$ are given by \eqref{eq1.5} with
$b(t)\in 2RV_{\beta, \rho_{2}}$, $c(t)\in 2RV_{\gamma, \rho_{3}}$,
$\beta, \gamma\in \R$, $\rho_{2}, \rho_{3}\leq 0$. These results are
more complex than results derived in Theorem \ref{th1}. In order to
state the main result Theorem \ref{th4} clearly, we first present
the following proposition.

\begin{proposition}\label{th3}
Let $F,G \in \mathcal{L}_{\alpha}$, $\alpha>0$.
Assume that $b(t)=e^{\alpha t}\overline{F}(t)\in 2RV_{\beta,\rho_{2}}$ with
auxiliary function $A_{2}(t)$ for $\rho_{2}\leq 0$, $\beta \in \R$,
and $c(t)=e^{\alpha t}\overline{G}(t)\in RV_{\gamma}$, $\gamma \in \R$. Then
\begin{itemize}
\item[(i)] if $\gamma\leq -1$, for large $t$ we have
\begin{eqnarray} \label{eq2.6}
\int_{0}^{\frac{t}{2}} \overline{F}(t-u)d G(u) =\overline{F}(t)M_{1}(\beta, \gamma, t)
\end{eqnarray}
with
\begin{eqnarray} \label{eq2.7}
M_{1}(\beta,\gamma, t)= \left\{
\begin{array}{ll}
\int_{0}^{\frac{t}{2}} c(u)du \left( \alpha
+ \frac{\alpha t c(t)}{\int_{0}^{\frac{t}{2}}c(u)du}
\left( \beta\int_{0}^{\frac{1}{2}} (1-u)^{\beta-1}(\log u)du
+(\log 2)\left( 1-2^{-\beta}\right) \right) \right. \\
\left.
+ \frac{1}{\int_{0}^{\frac{t}{2}}c(u)du}
+ o\left( \frac{1+t c(t)}{\int_{0}^{\frac{t}{2}}c(u)du} +A_{2}(t) \right)\right),
\quad \qquad \qquad \qquad \;\;\;
\gamma=-1 \, \text{with} \, m_{G}(\alpha)=\infty; \\
m_{G}(\alpha) -\alpha \int_{\frac{t}{2}}^{\infty} c(u)du
+ o\left( \int_{\frac{t}{2}}^{\infty} c(u)du +A_{2}(t) \right),
\qquad \quad \;\; \gamma=-1\, \text{with} \, m_{G}(\alpha)< \infty; \\
m_{G}(\alpha)+ \alpha\left( \int_{0}^{\frac{1}{2}}
\left( (1-u)^{\beta} -1 \right) u^{\gamma}du +\frac{1}{(1+\gamma)2^{1+\gamma}} \right) t c(t)+o(t c(t)+ A_{2}(t) ),\\
\qquad \qquad \qquad \qquad \qquad \qquad \qquad \qquad \qquad \qquad \qquad \qquad \qquad \qquad \quad \;
-2<\gamma<-1;\\
m_{G}(\alpha) - \alpha \beta t^{-1}\int_{1}^{\frac{t}{2}}u c(u)du +
o\left( t^{-1}\int_{1}^{\frac{t}{2}}u c(u)du + A_{2}(t) \right),\\
\qquad \qquad \qquad \qquad \qquad \qquad \qquad \qquad \qquad \qquad \qquad\;\;
\gamma=-2 \, \text{with} \, \int_{1}^{\infty} u c(u) du=\infty;\\
m_{G}(\alpha) - \beta t^{-1} \int_{0}^{\infty} ue^{\alpha u}d G(u)  +o(t^{-1}+ A_{2}(t)),
\qquad\;
\gamma \leq -2 \, \text{with} \, \int_{1}^{\infty} u c(u)du<\infty.
\end{array}
\right.
\end{eqnarray}

\item[(ii)] if $\gamma>-1$, suppose $c(t)\in 2RV_{\gamma,\rho_{3}}$
with auxiliary function $A_{3}(t)$ and $\rho_{3}\leq 0$,
$\gamma+\rho_{3}+1>0$, then \eqref{eq2.6} holds with
\begin{eqnarray} \label{eq2.11}
M_{1}(\beta, \gamma, t)
&=& t c(t)\left( \alpha\int_{0}^{\frac{1}{2}} (1-u)^{\beta}u^{\gamma}du
+ \frac{\alpha}{\rho_{2}}\int_{0}^{\frac{1}{2}} \left( (1-u)^{\rho_{2}} -1 \right) (1-u)^{\beta}u^{\gamma}du A_{2}(t) \right.
\nonumber\\
& &
\left.
+ \left( \beta\int_{0}^{\frac{1}{2}} (1-u)^{\beta-1}u^{\gamma}(2u-1)du
-2(1+\gamma)\int_{0}^{\frac{1}{2}}(1-u)^{\beta}u^{\gamma}du \right)t^{-1} \right. \nonumber\\
& &
\left.
+ \alpha\left( \frac{2^{-\rho_{3}}\beta }{\rho_{3}(\gamma+\rho_{3}+1)}
\int_{0}^{\frac{1}{2}} (1-u)^{\beta -1}u^{\gamma+1}\left((2u)^{\rho_{3}}-1\right)du  \right. \right.\nonumber\\
& &
\left.
\left.
+ \frac{2^{-\rho_{3}}(\gamma+1)-\gamma-\rho_{3}-1}
{\rho_{3}(\gamma+\rho_{3}+1)}
\int_{0}^{\frac{1}{2}}(1-u)^{\beta}u^{\gamma}du\right)A_{3}(t) + o\left( t^{-1}+A_{2}(t) +A_{3}(t) \right)\right)\nonumber\\
\end{eqnarray}
for large $t$.
\end{itemize}
\end{proposition}

\begin{theorem}\label{th4}
Let $F,G\in \mathcal{L}_{\alpha}$, $\alpha>0$. Assume that
$b(t)=e^{\alpha t}\overline{F}(t)\in 2RV_{\beta,\rho_{2}}$ with auxiliary function $A_{2}(t)$ for $\rho_{2}\leq 0$, $\beta \in \R$, and $c(t)=e^{\alpha t}\overline{G}(t)\in 2RV_{\gamma,\rho_{3}}$ with auxiliary function $A_{3}(t)$ for $\rho_{3}\leq 0$, $\gamma \in \R$. Then
\begin{equation}\label{eq2.16}
\overline{F*G}(t)=\overline{F}(t)M_{1}(\beta,\gamma,t) + \overline{G}(t)M_{2}(\beta,\gamma,t)
\end{equation}
for large $t$, where $M_{1}(\beta,\gamma, t)$ is given by Proposition \ref{th3},
and $M_{2}(\beta,\gamma, t)$ is given as follows:
\begin{itemize}

\item[(i)] if $\beta\leq -1$, for large $t$ we have
\begin{eqnarray}\label{eq2.17}
M_{2}(\beta,\gamma, t)=\left\{
\begin{array}{ll}
\int_{0}^{\frac{t}{2}} b(u)du \left( \alpha
+ \frac{\alpha t b(t)}{\int_{0}^{\frac{t}{2}}b(u)du}
\left( \gamma \int_{0}^{\frac{1}{2}} (1-u)^{\gamma-1}(\log u)du
+(\log 2)\left( 1-2^{-\gamma}\right) \right) \right. \\
\left.
+ \frac{1}{\int_{0}^{\frac{t}{2}}b(u)du}
+ o\left( \frac{1+t b(t)}{\int_{0}^{\frac{t}{2}}b(u)du} +A_{3}(t) \right)\right),
\qquad \qquad \qquad \;\;\;
\beta=-1 \, \text{with} \, m_{F}(\alpha)=\infty;\\

m_{F}(\alpha) -\alpha \int_{\frac{t}{2}}^{\infty} b(u)du
+ o\left( \int_{\frac{t}{2}}^{\infty} b(u)du +A_{3}(t) \right),
 \qquad  \;\;
\beta=-1\, \text{with} \, m_{F}(\alpha)< \infty; \\
m_{F}(\alpha)+ \alpha\left( \int_{0}^{\frac{1}{2}}
\left( (1-u)^{\gamma} -1 \right) u^{\beta}du +\frac{1}{(1+\beta)2^{1+\beta}} \right) t b(t)+o(t b(t)+A_{3}(t)),\\
\qquad \qquad \qquad \qquad \qquad \qquad \qquad \qquad \qquad \qquad \qquad \qquad  \qquad \qquad \;
-2<\beta<-1;\\
m_{F}(\alpha) - \alpha \gamma t^{-1} \int_{1}^{\frac{t}{2}}u b(u)du +
o\left( t^{-1}\int_{1}^{\frac{t}{2}}u b(u)du + A_{3}(t) \right),\\
\qquad \qquad \qquad \qquad \qquad \qquad \qquad\qquad \qquad\qquad \quad \;\;
\beta=-2 \, \text{with} \, \int_{1}^{\infty} u b(u) du=\infty;\\
m_{F}(\alpha) - \gamma t^{-1} \int_{0}^{\infty} ue^{\alpha u}d F(u)  +o(t^{-1}+A_{3}(t)), \quad \;
\beta \leq -2 \, \text{with} \, \int_{1}^{\infty} u b(u)du<\infty.
\end{array}
\right.
\end{eqnarray}

\item[(ii)] if $\beta>-1$ with $\beta+\rho_{2}+1>0$, for large $t$ we have
\begin{eqnarray}\label{eq2.18}
M_{2}(\beta,\gamma, t)
&=&
t b(t)\left( \alpha
\int_{0}^{\frac{1}{2}} (1-u)^{\gamma} u^{\beta}du
+ \frac{\alpha}{\rho_{3}}\int_{0}^{\frac{1}{2}} \left( (1-u)^{\rho_{3}} -1 \right) (1-u)^{\gamma}u^{\beta}du A_{3}(t)
\right. \nonumber\\
& &
\left.
+ \left( \gamma\int_{0}^{\frac{1}{2}} (1-u)^{\gamma-1}u^{\beta}(2u-1)du
-2(1+\beta)\int_{0}^{\frac{1}{2}}(1-u)^{\gamma}u^{\beta}du +2^{-\beta-\gamma}\right)t^{-1} \right. \nonumber\\
& &
\left.
+ \alpha\left( \frac{2^{-\rho_{2}}\gamma }{\rho_{2}(\beta+\rho_{2}+1)}
\int_{0}^{\frac{1}{2}} (1-u)^{\gamma -1}u^{\beta+1}\left((2u)^{\rho_{2}}-1\right)du  \right. \right.\nonumber\\
& &
\left.
\left.
+ \frac{2^{-\rho_{2}}(\beta+1)-\beta-\rho_{2}-1}
{\rho_{2}(\beta+\rho_{2}+1)}
\int_{0}^{\frac{1}{2}}(1-u)^{\gamma}u^{\beta}du\right)A_{2}(t)
 + o\left( t^{-1}+A_{2}(t) +A_{3}(t) \right)\right).\nonumber\\
\end{eqnarray}

\end{itemize}

\end{theorem}

\begin{remark}
For $\overline{F}(t)=b(t)e^{-\alpha t}, b(t)\in RV_{\beta}, \beta \in \R$,
it is clearly that $F\in \L_{\alpha}$ for all $\beta \in \R$. Note that
$F\in \S_{\alpha}$ if and only if $\beta \leq -1$ with $m_{F}(\alpha)<\infty$.
So Theorem \ref{th4} also derives the second-order tail asymptotics of convolution of distributions from a subclass of $\S_{\alpha}$.
\end{remark}

\begin{example}
Suppose $X$ and $Y$ are nonnegative random variables with distribution $F$
and $G$, respectively. Let $\overline{F}(t)=b(t)e^{-\alpha t}$ and
$\overline{G}(t)=c(t)e^{-\alpha t}$ with $b(t)=t^{-3}\left( 1+\frac{1}{4}A_{2}(t) +o(A_{2}(t)) \right)$, $|A_{2}(t)|\in RV_{-4}$,
$c(t)=\frac{\alpha^{\zeta-1}}{\Gamma(\zeta)}t^{\zeta-1}
\left(1+\frac{\zeta-1}{\alpha}t^{-1}+o(t^{-1}) \right)$,
 $\alpha>0$, $\zeta>1$. Here $\Gamma(\cdot)$
denotes the gamma function. From Lemma \ref{lem2.4}, it follows that
$b(t)\in 2RV_{-3, -4}$, $c(t)\in 2RV_{\zeta-1,-1}$ which imply that
$m_{F}(\alpha)<\infty$, $m_{G}(\alpha)=\infty$. Then second-order expansion
\eqref{eq2.16} leads to
\[\overline{F*G}(t)=\overline{G}(t)\left( m_{F}(\alpha)
-(\zeta-1)t^{-1}\int_{0}^{\infty} u e^{\alpha u}dF(u) +o(t^{-1}) \right),\]
\[\overline{F*F}(t)=2\overline{F}(t)\left( m_{F}(\alpha)
+ 3t^{-1}\int_{0}^{\infty} u e^{\alpha u}dF(u) +o(t^{-1}) \right) \]
and
\begin{eqnarray*}
\overline{G*G}(t)=
\frac{2\alpha^{\zeta-1}}{\Gamma(\zeta)}t^{\zeta}\overline{G}(t)
\left[ \alpha \int_{0}^{\frac{1}{2}} (1-y)^{\zeta-1}y^{\zeta-1} dy
+\left( (\zeta-1)\int_{0}^{\frac{1}{2}}(1-y)^{\zeta-2}y^{\zeta}dy
 + 2^{-2(\zeta-1)-1}\right)t^{-1} +o(t^{-1})  \right]
\end{eqnarray*}
for large $t$.
\end{example}

\section{Proofs  }
\label{sec4}

The aim of this section is to prove our main results. In order to
prove the main results, we give some auxiliary lemmas first. Without
loss of generality, we assume that the auxiliary functions of $2RV$
functions are positive eventually in the following proofs.

\begin{lemma}\label{lem2}
Let $\chi(t)\in 2RV_{\rho,\rho_{1}}$ with auxiliary function
$A_{1}(t)$ for $0<\rho<1$, $\rho_{1}\leq 0$. Assume that $\chi'(t)$
is nonincreasing eventually. Then for large $t$, we have
\begin{equation}\label{eq3.2}
\left| \frac{t\chi'(t)}{\chi(t)} -\rho \right|
\leq c \max\left( A_{1}^{c_{1}}(t), A_{1}^{1-c_{1}}(t) \right),
\end{equation}
where $c$ and $c_{1}$ are positive constants with $0<c_{1}<1$.
\end{lemma}

\noindent {\bf Proof.} First we consider the case of $\rho_{1}<0$.
Note that
\begin{equation}\label{eq3.3}
1+\rho x+\frac{\rho(\rho-1)x^{2}}{2}< (1+x)^{\rho}< 1+\rho x,
\end{equation}
\begin{equation}\label{eq3.4}
1-\rho x+ \left( 1-c_{2}^{-1} \right)^{\rho-2}
\frac{\rho(\rho-1)x^{2}}{2}< (1-x)^{\rho}< 1- \rho x,
\end{equation}
\begin{equation}\label{eq3.5}
1+\rho_{1} x< (1+x)^{\rho_{1}}< 1+\rho_{1} x+\frac{\rho_{1}(\rho_{1}-1)x^{2}}{2}
\end{equation}
and
\begin{equation}\label{eq3.6}
1-\rho_{1} x< (1-x)^{\rho_{1}}< 1-\rho_{1}x
+\left( 1-c_{2}^{-1} \right)^{\rho_{1}-2}\frac{\rho_{1}(\rho_{1}-1)x^{2}}{2}
\end{equation}
for $0<\rho<1$, $\rho_{1}<0$, $0<x<c_{2}^{-1}$, $1<c_{2}<1/\left(1-\left(2-2^{\rho}\right)^{\frac{1}{\rho}}\right)$.
By using \eqref{eq3.1}, \eqref{eq3.3}-\eqref{eq3.6}, we can get
\begin{eqnarray*}
\frac{\chi'(t)t}{\chi(t)}
&\leq& \frac{\int_{t\left( 1-A_{1}^{c_{1}}(t) \right)}^{t}\chi'(y)dy}{A_{1}^{c_{1}}(t)\chi(t)} \\
&\leq& A_{1}^{-c_{1}}(t)\left( 1- \left( 1-A_{1}^{c_{1}}(t) \right)^{\rho}
- A_{1}(t)\left( \left( 1-A_{1}^{c_{1}}(t) \right)^{\rho} \frac{\left( 1-A_{1}^{c_{1}}(t) \right)^{\rho_{1}}-1 }{\rho_{1}} -\varepsilon \left( 1-A_{1}^{c_{1}}(t) \right)^{\rho+\rho_{1}-\delta}  \right) \right) \\
&\leq& \rho-\left( 1-c_{2}^{-1} \right)^{\rho-2}\frac{\rho(\rho-1)}{2}A_{1}^{c_{1}}(t)
+ A_{1}(t)\left( 1-A_{1}^{c_{1}}(t) \right)^{\rho}\left( 1-\left( 1-c_{2}^{-1} \right)^{\rho_{1}-2}\frac{\rho_{1}-1}{2}A_{1}^{c_{1}}(t) \right)\\
& & +\varepsilon A_{1}^{1-c_{1}}(t)\left( 1-A_{1}^{c_{1}}(t)
\right)^{\rho+\rho_{1}-\delta}
\end{eqnarray*}
for large $t$. Similarly,
\begin{eqnarray*}
\frac{\chi'(t)t}{\chi(t)} &\geq&  \rho+
\frac{\rho(\rho-1)}{2}A_{1}^{c_{1}}(t) + A_{1}(t)\left(
1+A_{1}^{c_{1}}(t) \right)^{\rho}\left(
1+\frac{\rho_{1}-1}{2}A_{1}^{c_{1}}(t) \right) - \varepsilon
A_{1}^{1-c_{1}}(t)\left( 1+A_{1}^{c_{1}}(t)
\right)^{\rho+\rho_{1}+\delta}
\end{eqnarray*}
for large $t$, which implies
\begin{equation*}
\left| \frac{t\chi'(t)}{\chi(t)} -\rho \right|
\leq c \max\left( A_{1}^{c_{1}}(t), A_{1}^{1-c_{1}}(t) \right)
\end{equation*}
with $0<c_{1}<1$, $c>0$.

The arguments for the case of $\rho=0$ are similar. \qed

\begin{lemma}\label{lem3}
Let $\overline{F}(t)=e^{-\alpha t+\chi(t)}$ such that $\chi'(t)$ is
nonincreasing eventually. Suppose $\chi(t)\in 2RV_{\rho,\rho_{1}}$
with auxiliary function $A_{1}(t)$  and $0<\rho<1$,
$\rho+\rho_{1}<0$. Then
\begin{eqnarray}
\frac{\chi^{\frac{1}{2}}(t)\int_{f(2t)}^{t}\overline{F}(2t-u)d F(u)}
{t\overline{F}^{2}(t)}
&=& \frac{\alpha}{2}\sqrt{\frac{\pi}{\rho(1-\rho)}}
-\frac{\alpha\sqrt{\pi}(\rho-2)(\rho-3)}{32(\rho(1-\rho))^{\frac{3}{2}}}
\frac{1}{\chi(t)}\nonumber\\
& &
-\frac{1}{2}\sqrt{\frac{\rho\pi}{1-\rho}}\frac{\chi(t)}{t}
-\frac{\chi^{\frac{1}{2}}(t)}{2t}
+o\left( A_{1}(t)\chi(t)+\frac{1}{\chi(t)}+ \frac{\chi^{\frac{1}{2}}(t)}{t} \right)
\end{eqnarray}
for large $t$, where $f(t)=1/\chi'(t)$.
\end{lemma}

\noindent {\bf Proof.} From Lemma \ref{lem2.2}, we have
$\lim_{t\to \infty}t\chi'(t)/\chi(t)=\rho$, which implies that $\chi(t)$
is increasing eventually. For sufficiently large $t$ and $A$ satisfied
$f(2t)<A<t$, we can get
\begin{eqnarray}\label{eq3.8}
& &
\int_{\frac{t-A}{t}\chi^{\frac{1}{2}}(t)}^{\frac{t-f(2t)}{t}\chi^{\frac{1}{2}}(t)}
\left( \alpha-\chi'\left(  t-\frac{t y}{\chi^{\frac{1}{2}}(t)}  \right) \right)  \exp\left( \chi\left( t+\frac{t y}{\chi^{\frac{1}{2}}(t)}   \right)
+\chi\left( t-\frac{t y}{\chi^{\frac{1}{2}}(t)} \right) -2\chi(t) \right) dy \nonumber \\
&\leq&
\frac{\chi^{\frac{1}{2}}(t)}{t}\int_{0}^{A} \left( \alpha -\chi'(s) \right)
\exp\left( \chi\left( 2t-s\right) +\chi(s) -2\chi(t)\right) ds \nonumber\\
&\leq&
\left( \max_{s\in(0,A)}e^{\chi(s)} \right)\alpha A\frac{\chi^{\frac{1}{2}}(t)}{t}e^{\chi(2t)-2\chi(t)}
- \left(\min_{s\in(0,A)}\chi'(s)e^{\chi(s)} \right)A
\frac{\chi^{\frac{1}{2}}(t)}{t}e^{\chi(2t-A)-2\chi(t)} \nonumber\\
&\leq&
\left( \max_{s\in(0,A)}e^{\chi(s)} \right)\alpha A\frac{\chi^{\frac{1}{2}}(t)}{t}
e^{(1-\varepsilon)\left( 2^{\rho}-2 \right)\chi(t)}
- \left(\min_{s\in(0,A)}\chi'(s)e^{\chi(s)} \right)A
\frac{\chi^{\frac{1}{2}}(t)}{t}e^{(1+\varepsilon)\left( 2^{\rho}-2 \right)\chi(t)} \nonumber\\
&=&o\left( A_{1}(t)\chi(t)+\chi^{-1}(t)+ \frac{\chi^{\frac{1}{2}}(t)}{t} \right)
\end{eqnarray}
and
\begin{eqnarray}\label{eq3.9}
& &
\int^{\frac{t-A}{t}\chi^{\frac{1}{2}}(t)}_{\frac{1}{c_{2}}\chi^{\frac{1}{2}}(t)}
\left( \alpha-\chi'\left(  t-\frac{t y}{\chi^{\frac{1}{2}}(t)}  \right) \right)
\exp\left( \chi\left( t +\frac{t y}{\chi^{\frac{1}{2}}(t)}  \right)
+\chi\left( t-\frac{t y}{\chi^{\frac{1}{2}}(t)} \right) -2\chi(t) \right) dy \nonumber \\
&\leq&
\left(\alpha-\chi'\left( t(1-c_{2}^{-1}) \right) \right)
\left( 1-c_{2}^{-1}-At^{-1}\right)\chi^{\frac{1}{2}}(t)
\exp\Big( \chi\left( 2t-A \right)
+\chi\left( t(1-c_{2}^{-1}) \right)  -2\chi(t) \Big) \nonumber\\
&\leq&
\left(\alpha-\chi'\left( t(1-c_{2}^{-1}) \right) \right)
\left( 1-c_{2}^{-1}-At^{-1}\right)\chi^{\frac{1}{2}}(t)
\exp\Big( (1-\varepsilon)\left( 2^{\rho}+(1-c_{2}^{-1})^{\rho}-2 \Big)\chi(t) \right)\nonumber\\
&=&
o\left( A_{1}(t)\chi(t)+\chi^{-1}(t)+ \frac{\chi^{\frac{1}{2}}(t)}{t} \right),
\end{eqnarray}
where $1<c_{2}<1/\left( 1-(2-2^{\rho})^{\frac{1}{\rho}} \right)$.

In order to derived the desired result, we need some more precise inequalities as follows:
\begin{eqnarray}\label{eq3.10}
&&1+\rho x+\frac{\rho(\rho-1)}{2}x^{2}+\frac{\rho(\rho-1)(\rho-2)}{3!}x^{3}
+\frac{\rho(\rho-1)(\rho-2)(\rho-3)}{4!}x^{4} \nonumber\\
&<&(1+x)^{\rho} \nonumber\\
&<&1+\rho x+\frac{\rho(\rho-1)}{2}x^{2}+\frac{\rho(\rho-1)(\rho-2)}{3!}x^{3}
+\frac{\rho(\rho-1)(\rho-2)(\rho-3)}{4!}x^{4} \nonumber\\
& &+\frac{\rho(\rho-1)(\rho-2)(\rho-3)(\rho-4)}{5!}x^{5}
\end{eqnarray}
and
\begin{eqnarray}\label{eq3.11}
& &1-\rho x+\frac{\rho(\rho-1)}{2}x^{2}-\frac{\rho(\rho-1)(\rho-2)}{3!}x^{3}
+\frac{\rho(\rho-1)(\rho-2)(\rho-3)}{4!}x^{4} \nonumber\\
& &
-(1-c_{2}^{-1})^{\rho-5}\frac{\rho(\rho-1)(\rho-2)(\rho-3)(\rho-4)}{5!}x^{5} \nonumber\\
&<&(1-x)^{\rho} \nonumber\\
&<&1-\rho x+\frac{\rho(\rho-1)}{2}x^{2}-\frac{\rho(\rho-1)(\rho-2)}{3!}x^{3}
+\frac{\rho(\rho-1)(\rho-2)(\rho-3)}{4!}x^{4}
\end{eqnarray}
for $0<x<c_{2}^{-1}$.

Combining with \eqref{eq3.1}, \eqref{eq3.3}-\eqref{eq3.6}, we have
\begin{eqnarray}\label{eq3.12}
& &
\rho(\rho-1)y^{2}+\frac{\rho(\rho-1)(\rho-2)(\rho-3)}{12}\frac{y^{4}}{\chi(t)}
-(1-c_{2}^{-1})^{\rho-5}\frac{\rho(\rho-1)(\rho-2)(\rho-3)(\rho-4)}{5!}
\frac{y^{5}}{\chi^{\frac{3}{2}}(t)} \nonumber\\
& &
+\frac{1}{2\rho_{1}}\Big( (\rho+\rho_{1})(\rho+\rho_{1}-1)\left( 1+(1-c_{2}^{-1})^{\rho+\rho_{1}-2} \right) -\rho(\rho-1)\left(  1+(1-c_{2}^{-1})^{\rho-2} \right)\Big)y^{2}A_{1}(t) \nonumber\\
& &
-\varepsilon A_{1}(t)\chi(t)\left( \left( 1+\frac{y}{\chi^{\frac{1}{2}}(t)} \right)^{\rho+\rho_{1}+\varepsilon}
+\left( 1-\frac{y}{\chi^{\frac{1}{2}}(t)} \right)^{\rho+\rho_{1}-\varepsilon}  \right) \nonumber\\
&\leq&
\chi\left( t+\frac{yt}{\chi^{\frac{1}{2}}(t)}  \right)
+\chi\left( t-\frac{yt}{\chi^{\frac{1}{2}}(t)}  \right) -2\chi(t) \nonumber\\
&\leq&
\rho(\rho-1)y^{2}+\frac{\rho(\rho-1)(\rho-2)(\rho-3)}{12}\frac{y^{4}}{\chi(t)}
+\frac{\rho(\rho-1)(\rho-2)(\rho-3)(\rho-4)}{5!}
\frac{y^{5}}{\chi^{\frac{3}{2}}(t)} \nonumber\\
& &
+ \varepsilon A_{1}(t)\chi(t)\left( \left( 1+\frac{y}{\chi^{\frac{1}{2}}(t)} \right)^{\rho+\rho_{1}+\varepsilon}
+\left( 1-\frac{y}{\chi^{\frac{1}{2}}(t)} \right)^{\rho+\rho_{1}-\varepsilon}  \right)
\end{eqnarray}
for large $t$ and $0<y<c_{2}^{-1}\chi^{\frac{1}{2}}(t)$.

Note that $\rho(\rho-1)(\rho-2)(\rho-3)y^{4}/(12\chi(t))
+\rho(\rho-1)(\rho-2)(\rho-3)(\rho-4)y^{5}/\left(5!\chi^{\frac{3}{2}}(t) \right)<0$
for $0<\rho<1$, $0<y<c_{2}^{-1}\chi^{\frac{1}{2}}(t)$.
By using \eqref{eq3.12} and the following equality
\begin{equation}\label{eqadd1}
1+x<e^{x}<1+x+\frac{x^{2}}{2},\, x<0,
\end{equation}
we can get
\begin{eqnarray*}
& &
\int^{\frac{1}{c_{2}}\chi^{\frac{1}{2}}(t)}_{0}
\exp\left( \chi\left( t+\frac{yt}{\chi^{\frac{1}{2}}(t)}  \right)
+\chi\left( t-\frac{yt}{\chi^{\frac{1}{2}}(t)}  \right) -2\chi(t) \right) dy  \\
&\leq&
\exp\Big( \varepsilon \left( 1+(1-c_{2}^{-1})^{\rho+\rho_{1}-\varepsilon} \right) A_{1}(t)\chi(t) \Big) \\
& &
\times
\int_{0}^{\frac{1}{c_{2}}\chi^{\frac{1}{2}}(t)}
\Bigg( 1+ \frac{\rho(\rho-1)(\rho-2)(\rho-3)}{12}\frac{y^{4}}{\chi(t)}
+\frac{\rho(\rho-1)(\rho-2)(\rho-3)(\rho-4)}{5!}\frac{y^{5}}{\chi^{\frac{3}{2}}(t)}
\\
& &
+ \left(\frac{\rho(\rho-1)(\rho-2)(\rho-3)}{12} \right)^{2}\frac{y^{8}}{2\chi^{2}(t)}
\left( 1+\frac{\rho-4}{10}\frac{y}{\chi^{\frac{1}{2}}(t)} \right)^{2} \Bigg) \exp\left(\rho(\rho-1)y^{2} \right)dy \\
&=&
\int_{0}^{\infty}\exp\left( \rho(\rho-1)y^{2} \right)dy
+ \frac{\rho(\rho-1)(\rho-2)(\rho-3)}{12\chi(t)}\int_{0}^{\infty}y^{4}\exp\left( \rho(\rho-1)y^{2} \right)dy  +o\left(\chi^{-1}(t) + A_{1}(t)\chi(t) \right)
\end{eqnarray*}
and
\begin{eqnarray*}
& &
\int^{\frac{1}{c_{2}}\chi^{\frac{1}{2}}(t)}_{0}
\exp\left( \chi\left( t+\frac{t y}{\chi^{\frac{1}{2}}(t)} \right)
+\chi\left( t-\frac{t y}{\chi^{\frac{1}{2}}(t)}  \right) -2\chi(t) \right) dy  \\
&\geq&
\exp\Big( -\varepsilon \left( 1+(1-c_{2}^{-1})^{\rho+\rho_{1}-\varepsilon} \right) A_{1}(t)\chi(t) \Big)
\int_{0}^{\frac{1}{c_{2}}\chi^{\frac{1}{2}}(t)}
\Bigg( 1+ \frac{\rho(\rho-1)(\rho-2)(\rho-3)}{12}\frac{y^{4}}{\chi(t)}
\\
& &
- \left( 1-c_{2}^{-1} \right)^{\rho-5}\frac{\rho(\rho-1)(\rho-2)(\rho-3)(\rho-4)}{5!}\frac{y^{5}}{\chi^{\frac{3}{2}}(t)}
+ \frac{1}{2\rho_{1}}\left( (\rho+\rho_{1})(\rho+\rho_{1}-1)\left(  1+ \left(1-c_{2}^{-1} \right)^{\rho+\rho_{1}-2}\right) \right.  \\
&&
\left.-\rho(\rho-1)\left( 1+\left( 1-c_{2}^{-1} \right)^{\rho-2} \right) \right)y^{2}A_{1}(t) \Bigg) \exp\left(\rho(\rho-1)y^{2} \right)dy \\
&=&
\int_{0}^{\infty}\exp\left( \rho(\rho-1)y^{2} \right)dy
+ \frac{\rho(\rho-1)(\rho-2)(\rho-3)}{12\chi(t)}\int_{0}^{\infty}y^{4}\exp\left( \rho(\rho-1)y^{2} \right)dy  +o\left(\chi^{-1}(t) + A_{1}(t)\chi(t) \right)
\end{eqnarray*}
for large $t$, if $\rho+\rho_{1}<0$. Hence,
\begin{eqnarray}\label{eq3.13}
& &
\int^{\frac{1}{c_{2}}\chi^{\frac{1}{2}}(t)}_{0}
\exp\left( \chi\left( t+\frac{yt}{\chi^{\frac{1}{2}}(t)}  \right)
+\chi\left( t-\frac{yt}{\chi^{\frac{1}{2}}(t)}  \right) -2\chi(t) \right) dy \nonumber \\
&=&
\int_{0}^{\infty}\exp\left( \rho(\rho-1)y^{2} \right)dy
+ \frac{\rho(\rho-1)(\rho-2)(\rho-3)}{12\chi(t)}\int_{0}^{\infty}y^{4}\exp\left( \rho(\rho-1)y^{2} \right)dy  +o\left(\chi^{-1}(t) + A_{1}(t)\chi(t) \right)\nonumber\\
&=&
\frac{1}{2}\sqrt{\frac{\pi}{\rho(1-\rho)}}
-\frac{\sqrt{\pi}(\rho-2)(\rho-3)}{32\left(\rho(1-\rho) \right)^{\frac{3}{2}}}\frac{1}{\chi(t)}
+o\left(\chi^{-1}(t) + A_{1}(t)\chi(t) \right)
\end{eqnarray}
for large $t$.

Note that $\max\left( A_{1}^{c_{1}}(t), A_{1}^{1-c_{1}}(t) \right)=o\left( \chi^{-\frac{1}{2}}(t) \right)$ as $t\to \infty$, if $-\rho/(2\rho_{1})<c_{1}<1+\rho/(2\rho_{1})<1$ and $\rho+\rho_{1}<0$.
Due to \eqref{eq3.0}, \eqref{eq3.1}, \eqref{eq3.2}, \eqref{eq3.6} and \eqref{eq3.13},
we obtain
\begin{eqnarray*}
& & \left| \int_{0}^{\frac{1}{c_{2}}\chi^{\frac{1}{2}}(t)}
\left( \frac{t\left( 1-\frac{y}{\chi^{\frac{1}{2}}(t)}\right)
\chi'\left( t-\frac{yt}{\chi^{\frac{1}{2}}(t)} \right)}
{\chi\left( t -\frac{yt}{\chi^{\frac{1}{2}}(t)} \right)}
 -\rho \right)\frac{\chi\left( t-\frac{yt}{\chi^{\frac{1}{2}}(t)} \right)}
 {\left( 1-\frac{y}{\chi^{\frac{1}{2}}(t)}\right)\chi(t)} \right. \\
& &
\left.
\times
 \exp\left( \chi\left( t+\frac{yt}{\chi^{\frac{1}{2}}(t)}  \right)
+\chi\left( t-\frac{yt}{\chi^{\frac{1}{2}}(t)}  \right) -2\chi(t) \right)dy \right| \\
&\leq&
c\int_{0}^{\frac{1}{c_{2}}\chi^{\frac{1}{2}}(t)}
\left( A_{1}^{c_{1}}\left( t-\frac{yt}{\chi^{\frac{1}{2}}(t)} \right) + A_{1}^{1-c_{1}}\left( t-\frac{yt}{\chi^{\frac{1}{2}}(t)} \right)  \right)
\left( 1-\frac{y}{\chi^{\frac{1}{2}}(t)}\right)^{\rho-1}
\left( 1+\varepsilon\left( 1-\frac{y}{\chi^{\frac{1}{2}}(t)}\right)^{-\varepsilon} \right) \\
& &
\times
\exp\left( \chi\left( t+\frac{yt}{\chi^{\frac{1}{2}}(t)} \right)
+\chi\left( t-\frac{yt}{\chi^{\frac{1}{2}}(t)}  \right) -2\chi(t) \right)dy\\
&=&
o\left( \chi^{-\frac{1}{2}}(t) \right)
\end{eqnarray*}
with $-\rho/(2\rho_{1})<c_{1}<1+\rho/(2\rho_{1})<1$, and
\begin{eqnarray*}
& &
\left| \int_{0}^{\frac{1}{c_{2}}\chi^{\frac{1}{2}}(t)}
\left( \frac{\chi\left( t-\frac{yt}{\chi^{\frac{1}{2}}(t)}  \right)}{\left( 1-\frac{y}{\chi^{\frac{1}{2}}(t)} \right)\chi(t)}
-1+\frac{(\rho-1)y}{\chi^{\frac{1}{2}}(t)}\right)
\exp\left( \chi\left( t+\frac{yt}{\chi^{\frac{1}{2}}(t)}  \right)
+\chi\left( t-\frac{yt}{\chi^{\frac{1}{2}}(t)}  \right) -2\chi(t) \right)dy   \right| \\
&\leq&
\int_{0}^{\frac{1}{c_{2}}\chi^{\frac{1}{2}}(t)}
\Bigg[ \left( 1-c_{2}^{-1} \right)^{\rho-3}\frac{(\rho-1)(\rho-2)}{2}
\frac{y^{2}}{\chi(t)} + A_{1}(t)\left( 1-\frac{y}{\chi^{\frac{1}{2}}(t)} \right)^{\rho-1} \frac{\left( 1-\frac{y}{\chi^{\frac{1}{2}}(t)} \right)^{\rho_{1}}+1}{-\rho_{1}} \\
&&
+\varepsilon A_{1}(t) \left( 1-\frac{y}{\chi^{\frac{1}{2}}(t)} \right)^{\rho+\rho_{1}-1-\varepsilon}   \Bigg]
\exp\left( \chi\left( t+\frac{yt}{\chi^{\frac{1}{2}}(t)}  \right)
+\chi\left( t-\frac{yt}{\chi^{\frac{1}{2}}(t)}  \right) -2\chi(t) \right)dy \\
&=&
o\left( \chi^{-\frac{1}{2}}(t) \right)
\end{eqnarray*}
as $t\to \infty$, which implies that
\begin{eqnarray}\label{eq3.14}
& &
\int_{0}^{\frac{1}{c_{2}}\chi^{\frac{1}{2}}(t)}
\chi'\left( t-\frac{yt}{\chi^{\frac{1}{2}}(t)}  \right)
\exp\left( \chi\left( t +\frac{yt}{\chi^{\frac{1}{2}}(t)}  \right)
+\chi\left( t -\frac{yt}{\chi^{\frac{1}{2}}(t)}  \right) -2\chi(t) \right)dy \nonumber\\
&=&
\frac{\chi(t)}{t}\int_{0}^{\frac{1}{c_{2}}\chi^{\frac{1}{2}}(t)}
\frac{t\chi'\left( t-\frac{yt}{\chi^{\frac{1}{2}}(t)}  \right) }{\chi(t)}\exp\left( \chi\left( t +\frac{yt}{\chi^{\frac{1}{2}}(t)}  \right)
+\chi\left( t-\frac{yt}{\chi^{\frac{1}{2}}(t)}  \right) -2\chi(t) \right)dy \nonumber\\
&=&
\rho\int_{0}^{\infty}\exp\left( \rho(\rho-1)y^{2} \right)dy\frac{\chi(t)}{t}
+\rho(1-\rho)\int_{0}^{\infty}y\exp\left( \rho(\rho-1)y^{2} \right)dy\frac{\chi^{\frac{1}{2}}(t)}{t}
+ o\left( \frac{\chi^{\frac{1}{2}}(t)}{t} \right) \nonumber \\
&=&
\frac{1}{2}\sqrt{\frac{\rho\pi}{1-\rho}}\frac{\chi(t)}{t}
+\frac{\chi^{\frac{1}{2}}(t)}{2t}
+ o\left( \frac{\chi^{\frac{1}{2}}(t)}{t}\right)
\end{eqnarray}
for large $t$.

Combining \eqref{eq3.8}-\eqref{eq3.9} and
\eqref{eq3.13}-\eqref{eq3.14}, we have
\begin{eqnarray*}
& &
\frac{\chi^{\frac{1}{2}}(t)\int_{f(2t)}^{t}\overline{F}(2t-u)d F(u)}
{t\overline{F}^{2}(t)}\\
&=&
\int_{0}^{\frac{t-f(2t)}{t}\chi^{\frac{1}{2}}(t)}
\left( \alpha-\chi'\left(  t-\frac{yt}{\chi^{\frac{1}{2}}(t)}  \right) \right)  \exp\left( \chi\left( t +\frac{y t}{\chi^{\frac{1}{2}}(t)}  \right)
+\chi\left( t -\frac{yt}{\chi^{\frac{1}{2}}(t)}  \right) -2\chi(t) \right) dy\\
&=& \frac{\alpha}{2}\sqrt{\frac{\pi}{\rho(1-\rho)}}
-\frac{\alpha\sqrt{\pi}(\rho-2)(\rho-3)}{32(\rho(1-\rho))^{\frac{3}{2}}}
\frac{1}{\chi(t)}
-\frac{1}{2}\sqrt{\frac{\rho\pi}{1-\rho}}\frac{\chi(t)}{t}
-\frac{\chi^{\frac{1}{2}}(t)}{2t}
+o\left( A_{1}(t)\chi(t)+\frac{1}{\chi(t)}+ \frac{\chi^{\frac{1}{2}}(t)}{t} \right)
\end{eqnarray*}
for large $t$, which complete the proof.
\qed

\begin{lemma}\label{lem4}
Under the conditions of Lemma \ref{lem3}, we have
\begin{equation}\label{eq3.15}
\frac{\overline{F}^{2}\left( \frac{t}{2} \right)}
{\int_{0}^{\frac{t}{2}}\overline{F}(t-u)dF(u)}
= \frac{2^{2-\frac{\rho}{2}}\sqrt{\rho(1-\rho)}\chi^{\frac{1}{2}}(t)}
{\alpha \sqrt{\pi}t}(1+o(1))
\end{equation}
as $t\to \infty$.
\end{lemma}

\noindent {\bf Proof.}
From \eqref{eq3.1}, \eqref{eq3.10} and \eqref{eq3.11}, it follows that
\begin{eqnarray*}
& &
\rho(\rho-1)\frac{4\chi\left( \frac{t}{2} \right)}{t^{2}}u^{2}
+\frac{4\rho(\rho-1)(\rho-2)(\rho-3)}{3}\frac{\chi\left( \frac{t}{2} \right)}{t^{4}}u^{4}
-\frac{4\rho(\rho-1)(\rho-2)(\rho-3)(\rho-4)}{(1-c_{2}^{-1})^{5-\rho}15}
\frac{\chi\left( \frac{t}{2} \right)}{t^{5}}u^{5} \nonumber\\
& &
+ A_{1}\left(\frac{t}{2}\right)\chi\left(\frac{t}{2}\right)\left(
\left( 1+\frac{2u}{t} \right)^{\rho} \frac{\left( 1+\frac{2u}{t} \right)^{\rho_{1}}-1}{\rho_{1}}
+\left( 1-\frac{2u}{t} \right)^{\rho} \frac{\left( 1-\frac{2u}{t} \right)^{\rho_{1}}-1}{\rho_{1}} \right. \\
& &
\left.
-\varepsilon \left( 1+\frac{2u}{t} \right)^{\rho+\rho_{1}+\varepsilon}
-\varepsilon \left( 1-\frac{2u}{t} \right)^{\rho+\rho_{1}-\varepsilon}
\right)\\
&\leq&
\chi\left( \frac{t}{2} + u \right)
+\chi\left( \frac{t}{2} -u \right) -2\chi\left(\frac{t}{2}\right) \nonumber\\
&\leq&
\rho(\rho-1)\frac{4\chi\left( \frac{t}{2} \right)}{t^{2}}u^{2}
+\frac{4\rho(\rho-1)(\rho-2)(\rho-3)}{3}\frac{\chi\left( \frac{t}{2} \right)}{t^{4}}u^{4}
+ \frac{4\rho(\rho-1)(\rho-2)(\rho-3)(\rho-4)}{15}
\frac{\chi\left( \frac{t}{2} \right)}{t^{5}}u^{5} \nonumber\\
& &
+ A_{1}\left(\frac{t}{2}\right)\chi\left(\frac{t}{2}\right)\left(
\left( 1+\frac{2u}{t} \right)^{\rho} \frac{\left( 1+\frac{2u}{t} \right)^{\rho_{1}}-1}{\rho_{1}}
+\left( 1-\frac{2u}{t} \right)^{\rho} \frac{\left( 1-\frac{2u}{t} \right)^{\rho_{1}}-1}{\rho_{1}} \right. \\
& &
\left.
+\varepsilon \left( 1+\frac{2u}{t} \right)^{\rho+\rho_{1}+\varepsilon}
+\varepsilon \left( 1-\frac{2u}{t} \right)^{\rho+\rho_{1}-\varepsilon}
\right)
\end{eqnarray*}
for $0<u<t/(2c_{2})$, $1<c_{2}<1/\left( 1-\left( 2-2^{\rho} \right)^{\frac{1}{\rho}} \right)$.
Combining with \eqref{eqadd1}, we have
\begin{eqnarray*}
& &
\int_{0}^{\frac{t}{2c_{2}}} \exp\left( \chi\left( \frac{t}{2} + u \right)
+\chi\left( \frac{t}{2} -u \right) -2\chi\left(\frac{t}{2}\right)  \right)du \\
&\leq&
\exp\left( \left( \varepsilon\left( 1+ (1-c_{2}^{-1})^{\rho+\rho_{1}-\varepsilon} \right)
-\frac{2+(1+c_{2}^{-1})^{\rho}+(1-c_{2}^{-1})^{\rho+\rho_{1}}}{\rho_{1}} \right)
A_{1}\left( \frac{t}{2} \right)\chi\left( \frac{t}{2} \right) \right) \\
& &
\times
\frac{t\chi^{-\frac{1}{2}}\left( \frac{t}{2} \right)}{2\sqrt{2\rho(1-\rho)}}
\int_{0}^{\frac{\chi^{\frac{1}{2}}\left( \frac{t}{2} \right)}{c_{2}(2\rho(1-\rho))^{-\frac{1}{2}}}}
\exp\left( -\frac{s^{2}}{2} -\frac{(\rho-2)(\rho-3)}
{48\rho(1-\rho)\chi\left(\frac{t}{2} \right)}  s^{4}
-\frac{(\rho-2)(\rho-3)(\rho-4)}
{480\sqrt{2}(\rho(1-\rho))^{\frac{3}{2}}\chi^{\frac{3}{2}}(t)}s^{5}
\right)ds\\
&\leq&
\frac{\left( 1+O\left( A_{1}\left( \frac{t}{2} \right)
\chi\left( \frac{t}{2} \right) \right) \right)t}{2\sqrt{2\rho(1-\rho)}\chi^{\frac{1}{2}}\left( \frac{t}{2} \right)}
\int_{0}^{\frac{\chi^{\frac{1}{2}}\left( \frac{t}{2} \right)}{c_{2}(2\rho(1-\rho))^{-\frac{1}{2}}}}
\left( 1-\frac{(\rho-2)(\rho-3)}
{48\rho(1-\rho)\chi\left(\frac{t}{2} \right)}  s^{4}
-\frac{(\rho-2)(\rho-3)(\rho-4)}
{480\sqrt{2}(\rho(1-\rho))^{\frac{3}{2}}\chi^{\frac{3}{2}}(t)}s^{5}
\right. \\
& &
\left.
+ \frac{1}{2}\left( \frac{(\rho-2)(\rho-3)}
{48\rho(1-\rho)\chi\left(\frac{t}{2} \right)} \right)^{2}s^{8}
\left( 1+\frac{\rho-4}{10\sqrt{2\rho(1-\rho)}\chi^{\frac{1}{2}}(t)}s \right)^{2}  \right) \exp\left( -\frac{s^{2}}{2} \right)ds \\
&=&
\frac{\sqrt{\pi}t}{4\sqrt{\rho(1-\rho)}\chi^{\frac{1}{2}}\left( \frac{t}{2} \right)}(1+o(1))
\end{eqnarray*}
and
\begin{eqnarray*}
& &
\int_{0}^{\frac{t}{2c_{2}}} \exp\left( \chi\left( \frac{t}{2} + u \right)
+\chi\left( \frac{t}{2} -u \right) -2\chi\left(\frac{t}{2}\right)  \right)du \\
&\geq&
\exp\left( -\left( \varepsilon\left( 1+ (1-c_{2}^{-1})^{\rho+\rho_{1}-\varepsilon} \right)
-\frac{2+(1+c_{2}^{-1})^{\rho}+(1-c_{2}^{-1})^{\rho+\rho_{1}}}{\rho_{1}} \right)
A_{1}\left( \frac{t}{2} \right)\chi\left( \frac{t}{2} \right) \right) \\
& &
\times
\frac{t\chi^{-\frac{1}{2}}\left( \frac{t}{2} \right)}{2\sqrt{2\rho(1-\rho)}}
\int_{0}^{\frac{\chi^{\frac{1}{2}}\left( \frac{t}{2} \right)}{c_{2}(2\rho(1-\rho))^{-\frac{1}{2}}}}
\exp\left( -\frac{s^{2}}{2} -\frac{(\rho-2)(\rho-3)}
{48\rho(1-\rho)\chi\left(\frac{t}{2} \right)}  s^{4}
+\frac{\left( 1-c_{2}^{-1} \right)^{\rho-5}(\rho-2)(\rho-3)(\rho-4)}
{480\sqrt{2}(\rho(1-\rho))^{\frac{3}{2}}\chi^{\frac{3}{2}}(t)}s^{5}
\right)ds\\
&=&
\frac{\sqrt{\pi}t}{4\sqrt{\rho(1-\rho)}\chi^{\frac{1}{2}}\left( \frac{t}{2} \right)}(1+o(1))
\end{eqnarray*}
for large $t$, if $0<\rho<1$ and $\rho+\rho_{1}<0$. So
\begin{equation}\label{eq3.16}
\int_{0}^{\frac{t}{2c_{2}}} \exp\left( \chi\left( \frac{t}{2} + u \right)
+\chi\left( \frac{t}{2} -u \right) -2\chi\left(\frac{t}{2}\right)  \right)du =\frac{\sqrt{\pi}t}{4\sqrt{\rho(1-\rho)}\chi^{\frac{1}{2}}\left( \frac{t}{2} \right)}(1+o(1))
\end{equation}
for large $t$.

By arguments similar to \eqref{eq3.8}, \eqref{eq3.9},
\eqref{eq3.14}, we can get
\begin{equation*}
\int_{0}^{\frac{t}{2c_{2}}} \chi'\left( \frac{t}{2} -u \right)\exp\left( \chi\left( \frac{t}{2} + u \right)
+\chi\left( \frac{t}{2} -u \right) -2\chi\left(\frac{t}{2}\right)  \right)du =\frac{1}{2}\sqrt{\frac{\rho\pi}{1-\rho}}\chi^{\frac{1}{2}}\left(\frac{t}{2}\right)(1+o(1))
\end{equation*}
and
\begin{equation*}
\int^{\frac{t}{2}}_{\frac{t}{2c_{2}}} \left( \alpha-\chi'\left( \frac{t}{2} -u \right) \right) \exp\left( \chi\left( \frac{t}{2} + u \right)
+\chi\left( \frac{t}{2} -u \right) -2\chi\left(\frac{t}{2}\right)  \right)du =o\left( \frac{t}{\chi^{\frac{1}{2}}(t)} \right)
\end{equation*}
for large $t$. Combining with \eqref{eq3.16}, we have
\begin{eqnarray*}
\frac{\int_{0}^{\frac{t}{2}}\overline{F}(t-u)dF(u)}
{\overline{F}^{2}\left( \frac{t}{2} \right)}
&=&
\frac{\alpha\sqrt{\pi}t}{4\sqrt{\rho(1-\rho)}\chi^{\frac{1}{2}}\left( \frac{t}{2} \right)}(1+o(1))\\
&=& \frac{2^{\frac{\rho}{2}-2}\alpha \sqrt{\pi}t}
{\sqrt{\rho(1-\rho)}\chi^{\frac{1}{2}}(t)}(1+o(1))
\end{eqnarray*}
for large $t$, which deduces the desired result.

The proof is complete.
\qed

\begin{lemma}\label{lem5}
Assume that $b(t)\in 2RV_{\beta, \rho_{2}}$ with auxiliary function
$A_{2}(t)$ and $\rho_{2}\leq 0$, $\beta\in \R$, and $c(t)\in
2RV_{\gamma, \rho_{3}}$ with auxiliary function $A_{3}(t)$ and
$\rho_{3}\leq 0$, $\gamma\in \R$. Then $b(t)c(t)\in
2RV_{\beta+\gamma, \max(\rho_{2}, \rho_{3})}$.
\end{lemma}

\noindent {\bf Proof.} From Lemma \ref{lem2.3}, it follows that
\begin{eqnarray}\label{eq3.18}
\frac{b(tx)c(tx)}{b(t)c(t)}=x^{\beta+\gamma}\left( 1+A_{2}(t)\frac{x^{\rho_{2}}-1}{\rho_{2}} +A_{3}(t)\frac{x^{\rho_{3}}-1}{\rho_{3}}
+o\left(A_{2}(t)+A_{3}(t) \right) \right)
\end{eqnarray}
for large $t$ and fixed $x>0$, which implies the desired result.
\qed

\begin{lemma}\label{lem6}
Let $G\in \mathcal{L}_{\alpha}$ and
$\mathcal{C}(t)=\int_{0}^{t}e^{\alpha u}dG(u)$ with $\alpha>0$.
Assume that $c(t)=e^{\alpha t}\overline{G}(t)\in 2RV_{\gamma, \rho_{3}}$
with auxiliary function $A_{3}(t)$ for $\rho_{3}\leq 0$, $\gamma>-1$, $\gamma+\rho_{3}+1>0$,
then $\mathcal{C}(t) \in 2RV_{\gamma+1, \max(\rho_{3}, -1)}$.
\end{lemma}

\noindent {\bf Proof.}
By using Lemma \ref{lem2.3}, for every $\varepsilon>0$ and $x, y>0$, there exists
$t_{0}=t_{0}(\varepsilon)$ such that all $tx, txy\geq t_{0}$,
\begin{eqnarray*}
& &
\left| \frac{\int_{\frac{t_{0}}{tx}}^{1} \left( \frac{c(txy)}{c(tx)} -y^{\gamma} \right)dy}{A_{3}(tx)}
- \int_{\frac{t_{0}}{tx}}^{1}y^{\gamma}\frac{y^{\rho_{3}}-1}{\rho_{3}}dy \right| \\
&\leq&
\int_{\frac{t_{0}}{tx}}^{1}\left| \frac{\frac{c(txy)}{c(tx)} -y^{\gamma} }{A_{3}(tx)} - y^{\gamma}\frac{y^{\rho_{3}}-1}{\rho_{3}} \right|dy \\
&\leq&
\varepsilon \int_{\frac{t_{0}}{t x}}^{1}y^{\gamma+\rho_{3}-\varepsilon}dy\\
&=&
\frac{\varepsilon}{\gamma+\rho_{3}+1-\varepsilon}
\left( 1-\left(\frac{t_{0}}{t x}\right)^{\gamma+\rho_{3}+1-\varepsilon} \right),
\end{eqnarray*}
which implies
\begin{eqnarray*}
& &
\frac{1}{\gamma+1} -\frac{A_{3}(tx)}{(\gamma+1)(\gamma+\rho_{3}+1)}
- A_{3}(tx)\left( \frac{\varepsilon}{\gamma+\rho_{3}+1-\varepsilon}
-\frac{1}{\rho_{3}(\gamma+1)}\left( \frac{t_{0}}{tx} \right)^{\gamma+1} \right)
-\frac{1}{\gamma+1}\left( \frac{t_{0}}{tx} \right)^{\gamma+1}\\
&\leq&
\frac{\int_{0}^{tx} c(y)dy}{tx c(tx)}\\
&\leq&
\frac{1}{\gamma+1} -\frac{A_{3}(tx)}{(\gamma+1)(\gamma+\rho_{3}+1)}
+ A_{3}(tx)\left( \frac{\varepsilon}{\gamma+\rho_{3}+1-\varepsilon}
-\frac{1}{\rho_{3}(\gamma+\rho_{3}+1)}\left( \frac{t_{0}}{t x} \right)^{\gamma+\rho_{3}+1} \right)
+\frac{\int_{0}^{t_{0}}c(y)dy}{tx c(tx)}
\end{eqnarray*}
for large $t$.

Note that $\mathcal{C}(t)=1-c(t)+\alpha \int_{0}^{t}c(u)du$. Combining with Lemma \ref{lem2.3}, for large $t$ we have
\begin{eqnarray}\label{eq3.19}
& &
-\frac{\gamma+1}{\alpha t}x^{\gamma}
\left( 1+A_{3}(t)\left( \frac{x^{\rho_{3}}-1}{\rho_{3}} + \varepsilon x^{\rho_{3}} \max\left( x^{\varepsilon}, x^{-\varepsilon} \right) \right) \right)
+
\frac{\gamma+1}{\alpha t c(t)} \nonumber \\
& &
+ x^{\gamma+1}\left( 1+ A_{3}(t)\left( \frac{x^{\rho_{3}}-1}{\rho_{3}} - \varepsilon x^{\rho_{3}}\max\left( x^{\varepsilon}, x^{-\varepsilon} \right) \right) \right)
\times
\left[ 1 -\left( \frac{t_{0}}{tx} \right)^{\gamma+1} \right. \nonumber\\
& &
 \left. \;\; -A_{3}(t)x^{\rho_{3}}\left( 1+\varepsilon\max\left( x^{\varepsilon}, x^{-\varepsilon} \right) \right)\left( \frac{1}{\gamma+\rho_{3}+1}
+ \frac{\varepsilon(\gamma+1)}{\gamma+\rho_{3}+1-\varepsilon}
-\frac{1}{\rho_{3}}\left( \frac{t_{0}}{tx} \right)^{\gamma+1}\right)
 \right]
 \nonumber\\
&\leq&
\frac{\mathcal{C}(tx)}{\frac{\alpha}{\gamma+1}t c(t)}
= \frac{\gamma+1}{\alpha t c(t)} -\frac{\gamma+1}{\alpha t}\frac{c(tx)}{c(t)}
+ \frac{(\gamma+1)\int_{0}^{tx} c(u)du}{tx c(tx)}\frac{x c(tx)}{c(t)}\nonumber\\
&\leq&
-\frac{\gamma+1}{\alpha t}x^{\gamma}
\left( 1+A_{3}(t)\left( \frac{x^{\rho_{3}}-1}{\rho_{3}} - \varepsilon x^{\rho_{3}} \max\left( x^{\varepsilon}, x^{-\varepsilon} \right) \right) \right)
+\frac{(\gamma+1)\left( \alpha\int_{0}^{t_{0}}c(y)dy+1  \right)}
{\alpha t c(t)}
\nonumber\\
& &
+ x^{\gamma+1}\left( 1+ A_{3}(t)\left( \frac{x^{\rho_{3}}-1}{\rho_{3}} + \varepsilon x^{\rho_{3}}\max\left( x^{\varepsilon}, x^{-\varepsilon} \right) \right) \right)
\times
\left[ 1 - \frac{A_{3}(t)x^{\rho_{3}}}{\gamma+\rho_{3}+1}\left( 1-\varepsilon \max\left( x^{\varepsilon}, x^{-\varepsilon} \right) \right)
\right.  \nonumber\\
& &
\left. \;\;
+ A_{3}(t)x^{\rho_{3}}\left(  1+ \varepsilon \max\left( x^{\varepsilon}, x^{-\varepsilon} \right) \right)\left( \frac{\varepsilon (\gamma+1)}{\gamma+\rho_{3}+1-\varepsilon}
-\frac{\gamma+1}{\rho_{3}(\gamma+\rho_{3}+1)}\left( \frac{t_{0}}{tx} \right)^{\gamma+\rho_{3}+1} \right)
 \right]
.
\end{eqnarray}
Hence,
\begin{equation}\label{eq3.20}
\frac{\mathcal{C}(tx)}{\frac{\alpha}{\gamma+1}t c(t)}
= x^{\gamma+1}+A_{3}(t)x^{\gamma+1}\frac{x^{\rho_{3}}-1}{\rho_{3}}
-A_{3}(t)\frac{x^{\gamma+\rho_{3}+1}}{\gamma+\rho_{3}+1}
-\frac{\gamma+1}{\alpha t}x^{\gamma} + o\left( \frac{1}{t} + A_{3}(t) \right)
\end{equation}
for large $t$, if $\gamma>-1$, $\gamma+\rho_{3}+1>0$.
Then we can get
\begin{eqnarray*}
\frac{\mathcal{C}(tx)}{\mathcal{C}(t)}
&=&\frac{x^{\gamma+1}+A_{3}(t)x^{\gamma+1}\frac{x^{\rho_{3}}-1}{\rho_{3}}
-A_{3}(t)\frac{x^{\gamma+\rho_{3}+1}}{\gamma+\rho_{3}+1}
-\frac{\gamma+1}{\alpha t}x^{\gamma}
+ o\left( \frac{1}{t} + A_{3}(t) \right)}
{1-\frac{A_{3}(t)}{\gamma+\rho_{3}+1}-\frac{\gamma+1}{\alpha t}
+ o\left( \frac{1}{t} + A_{3}(t) \right)} \\
&=&
x^{\gamma+1}+ x^{\gamma+1}\left( x^{\rho_{3}}-1 \right)\frac{\gamma+1}{\rho_{3}(\gamma+\rho_{3}+1)}A_{3}(t)
+x^{\gamma}\left( x-1 \right) \frac{\gamma+1}{\alpha t}
+o\left( \frac{1}{t} + A_{3}(t) \right)
\end{eqnarray*}
for large $t$ and fixed $x>0$, which deduces that $\mathcal{C}(t)\in 2RV_{\gamma+1, \max(\rho_{3},-1)}$.

The proof is complete.
\qed

\begin{lemma}\label{lem7}
Let $F, G \in \mathcal{L}_{\alpha}$, $\alpha>0$. Assume that
$b(t)=e^{\alpha t}\overline{F}(t)\in 2RV_{\beta, \rho_{2}}$ with
auxiliary function $A_{2}(t)$ and $\rho_{2}\leq 0$, $\beta\in \R$,
and $c(t)=e^{\alpha t}\overline{G}(t)\in 2RV_{\gamma, \rho_{3}}$
with auxiliary function $A_{3}(t)$ and $\rho_{3}\leq 0$,
$\gamma>-1$, $\gamma +\rho_{3}+1>0$. Then
\begin{eqnarray*}
\frac{\int_{0}^{\frac{t}{2}} \overline{F}(t-u)d G(u)}{\overline{F}(t)t c(t)}
&=&
\alpha\int_{0}^{\frac{1}{2}} (1-u)^{\beta}u^{\gamma}du
+ \frac{\alpha}{\rho_{2}}\int_{0}^{\frac{1}{2}} \left( (1-u)^{\rho_{2}} -1 \right) (1-u)^{\beta}u^{\gamma}du A_{2}(t)
\nonumber\\
& &
+ \left( \beta\int_{0}^{\frac{1}{2}} (1-u)^{\beta-1}u^{\gamma}(2u-1)du
-2(1+\gamma)\int_{0}^{\frac{1}{2}}(1-u)^{\beta}u^{\gamma}du \right)t^{-1} \nonumber\\
& &
+ \alpha\left( \frac{2^{-\rho_{3}}\beta }{\rho_{3}(\gamma+\rho_{3}+1)}
\int_{0}^{\frac{1}{2}} (1-u)^{\beta -1}u^{\gamma+1}\left((2u)^{\rho_{3}}-1\right)du  \right. \nonumber\\
& &
\left.
+ \frac{2^{-\rho_{3}}(\gamma+1)-\gamma-\rho_{3}-1}
{\rho_{3}(\gamma+\rho_{3}+1)}
\int_{0}^{\frac{1}{2}}(1-u)^{\beta}u^{\gamma}du\right)A_{3}(t)
+ o\left( t^{-1}+A_{2}(t) +A_{3}(t) \right)
\end{eqnarray*}
for large $t$.
\end{lemma}

\noindent {\bf Proof.} Note that $\C(t)=\int_{0}^{t}e^{\alpha u}d
G(u)$. From \eqref{eq3.19}, \eqref{eq3.20} and the dominated
convergence theorem, we have
\begin{eqnarray*}
& &
\int_{\frac{1}{t}}^{\frac{1}{2}} (1-u)^{\beta-1}\left( \frac{\C(ut)}{\C\left( \frac{t}{2}\right)} -(2u)^{\gamma+1} \right)du \\
&=&
\frac{\frac{\alpha}{\gamma+1}\frac{t}{2}c\left(\frac{t}{2}\right)}
{\C\left( \frac{t}{2} \right)}
\int_{\frac{1}{t}}^{\frac{t}{2}} (1-u)^{\beta-1}
\left( \frac{\C\left( ut \right)}{\frac{\alpha}{\gamma+1}\frac{t}{2}c\left(\frac{t}{2}\right)}
-(2u)^{\gamma+1} \right)du
+ \left( \frac{\frac{\alpha}{\gamma+1}\frac{t}{2}c\left(\frac{t}{2}\right)}
{\C\left( \frac{t}{2} \right)}-1 \right)
\int_{\frac{1}{t}}^{\frac{t}{2}} (1-u)^{\beta-1}(2u)^{\gamma+1} du \\
&=&
\left( 1+\frac{A_{3}\left( \frac{t}{2} \right)}{\gamma+\rho_{3}+1}
+\frac{2(\gamma+1)}{\alpha t} + o\left( \frac{1}{t} +A_{3}(t) \right) \right)
\left[
-\frac{2(\gamma+1)}{\alpha t} \int_{0}^{\frac{1}{2}} (1-u)^{\beta-1}(2u)^{\gamma}du \right. \\
& &
\left.
+ A_{3}\left( \frac{t}{2} \right)\left( \int_{0}^{\frac{1}{2}} (1-u)^{\beta-1} (2u)^{\gamma+1}\frac{(2u)^{\rho_{3}}-1}{\rho_{3}}du
- \int_{0}^{\frac{1}{2}} (1-u)^{\beta-1} \frac{(2u)^{\gamma+\rho_{3}+1}}{\gamma+\rho_{3}+1}du \right) + o\left( \frac{1}{t} +A_{3}(t) \right)\right]\\
& &
+ \int_{0}^{\frac{1}{2}} (1-u)^{\beta-1}(2u)^{\gamma+1}du
\left( \frac{A_{3}\left( \frac{t}{2} \right)}{\gamma+\rho_{3}+1}
+\frac{2(\gamma+1)}{\alpha t} + o\left( \frac{1}{t} +A_{3}(t) \right) \right)\\
&=&
A_{3}(t)\frac{2^{\gamma-\rho_{3}+1}(\gamma+1)}{\rho_{3}(\gamma+\rho_{3}+1)}
\int_{0}^{\frac{1}{2}} (1-u)^{\beta-1}u^{\gamma+1}\left( (2u)^{\rho_{3}} -1 \right)du \\
& &
+\frac{2^{\gamma+1}(\gamma+1)}{\alpha t} \int_{0}^{\frac{1}{2}} (1-u)^{\beta-1}u^{\gamma}(2u-1)du + o\left( \frac{1}{t} +A_{3}(t) \right)
\end{eqnarray*}
for large $t$. So
\begin{eqnarray}\label{eq3.21}
& &
\frac{\int_{0}^{\frac{t}{2}}\left( 1-\frac{u}{t} \right)^{\beta}
e^{\alpha u}d G(u)}{\C\left( \frac{t}{2} \right)}
-(\gamma+1)2^{\gamma+1}\int_{0}^{\frac{1}{2}} (1-u)^{\beta}u^{\gamma}du \nonumber\\
&=&
\frac{\int_{0}^{1}\left( 1-\frac{u}{t} \right)^{\beta}
e^{\alpha u}d G(u)-\left( 1-\frac{1}{t}\right)^{\beta}\C(1)}
{\C\left( \frac{t}{2}\right)} -\beta\int_{0}^{\frac{1}{t}} (1-u)^{\beta-1}(2u)^{\gamma+1}du \nonumber\\
& &
+ \beta \int_{\frac{1}{t}}^{\frac{1}{2}} (1-u)^{\beta-1}\left( \frac{\C(ut)}{\C\left( \frac{t}{2}\right)} -(2u)^{\gamma+1} \right)du \nonumber\\
&=&
A_{3}(t)\frac{2^{\gamma-\rho_{3}+1}\beta(\gamma+1)}{\rho_{3}(\gamma+\rho_{3}+1)}
\int_{0}^{\frac{1}{2}} (1-u)^{\beta-1}u^{\gamma+1}\left( (2u)^{\rho_{3}} -1 \right)du  \nonumber \\
& &
+\frac{2^{\gamma+1}\beta(\gamma+1)}{\alpha t} \int_{0}^{\frac{1}{2}} (1-u)^{\beta-1}u^{\gamma}(2u-1)du + o\left( \frac{1}{t} +A_{3}(t) \right),
\end{eqnarray}
due to $\C(t)\in RV_{\gamma+1}$ and $\gamma+\rho_{3}+1>0$.

By using integration by parts, for $\gamma>-1$ we can get
\begin{equation}\label{eq3.22}
\frac{\int_{0}^{\frac{t}{2}}\left( 1-\frac{u}{t} \right)^{\beta}
\frac{\left( 1-\frac{u}{t} \right)^{\rho_{2}}-1}{\rho_{2}}
e^{\alpha u}d G(u)}{\C\left( \frac{t}{2} \right)}
\to \frac{\gamma+1}{\rho_{2}}2^{\gamma+1}\int_{0}^{\frac{1}{2}}
\left( (1-u)^{\rho_{2}} -1 \right)(1-u)^{\beta}u^{\gamma}du
\end{equation}
and
\begin{equation}\label{eq3.23}
\frac{\int_{0}^{\frac{t}{2}}\left( 1-\frac{u}{t} \right)^{\beta+\rho_{2}-\varepsilon}
e^{\alpha u}d G(u)}{\C\left( \frac{t}{2} \right)}
\to (\gamma+1)2^{\gamma+1}\int_{0}^{\frac{1}{2}}
(1-u)^{\beta+\rho_{2}-\varepsilon}u^{\gamma}du
\end{equation}
as $t\to \infty$.

Note that $\int_{0}^{\frac{t}{2}} \frac{\overline{F}(t-u)}
{\overline{F}(t)}d
G(u)=\int_{0}^{\frac{t}{2}}\frac{b(t-u)}{b(t)}e^{\alpha u}d G(u)$.
From Lemma \ref{lem2.3} and \eqref{eq3.20}-\eqref{eq3.23}, it
follows that
\begin{eqnarray*}
&&
\frac{\int_{0}^{\frac{t}{2}} \overline{F}(t-u)d G(u)}{\overline{F}(t)t c(t)}\\
&=&
\alpha\int_{0}^{\frac{1}{2}} (1-u)^{\beta}u^{\gamma}du
+ \left( \beta\int_{0}^{\frac{1}{2}} (1-u)^{\beta-1}u^{\gamma}(2u-1)du
-2(1+\gamma)\int_{0}^{\frac{1}{2}}(1-u)^{\beta}u^{\gamma}du \right)\frac{1}{t} \nonumber\\
& &
+ \alpha\left( \frac{\rho^{-1}_{3}2^{-\rho_{3}}\beta }{\gamma+\rho_{3}+1}
\int_{0}^{\frac{1}{2}} (1-u)^{\beta -1}u^{\gamma+1}\left((2u)^{\rho_{3}}-1\right)du + \frac{2^{-\rho_{3}}(\gamma+1)-\gamma-\rho_{3}-1}
{\rho_{3}(\gamma+\rho_{3}+1)}
\int_{0}^{\frac{1}{2}}(1-u)^{\beta}u^{\gamma}du\right)A_{3}(t) \nonumber\\
& &
+ \frac{\alpha}{\rho_{2}}\int_{0}^{\frac{1}{2}} \left( (1-u)^{\rho_{2}} -1 \right) (1-u)^{\beta}u^{\gamma}du A_{2}(t) + o\left( t^{-1}+A_{2}(t) +A_{3}(t) \right)
\end{eqnarray*}
for large $t$, which complete the proof.
\qed

\begin{lemma}\label{lem8}
Let $F, G\in \L_{\alpha}$, $\alpha>0$. Assume that $b(t)=e^{\alpha t}\overline{F}(t)\in 2RV_{\beta, \rho_{2}}$ with auxiliary function $A_{2}(t)$ for $\rho_{2}\leq 0$, $\beta\in R$, and $c(t)=e^{\alpha t}\overline{G}(t)\in RV_{\gamma}$ for $\gamma \leq -1$. For large $t$,
we have
\begin{itemize}
\item[(i)] if $\gamma=-1$ with $m_{G}(\alpha)=\infty$,
\begin{eqnarray}\label{eq3.24}
\frac{\int_{0}^{\frac{t}{2}}\overline{F}(t-u)d G(u)}{\overline{F}(t)\int_{0}^{\frac{t}{2}}c(u)du}
&=&\alpha
+ \frac{\alpha t c(t)}{\int_{0}^{\frac{t}{2}}c(u)du}
\left( \beta\int_{0}^{\frac{1}{2}} (1-u)^{\beta-1}(\log u)du
+(\log 2)\left( 1-2^{-\beta}\right) \right) \nonumber \\
& &
+ \frac{1}{\int_{0}^{\frac{t}{2}}c(u)du}
+ o\left( \frac{1+t c(t)}{\int_{0}^{\frac{t}{2}}c(u)du} +A_{2}(t) \right);
\end{eqnarray}

\item[(ii)] if $\gamma=-1$ with $m_{G}(\alpha)<\infty$,
\begin{eqnarray}\label{eq3.25}
\frac{\int_{0}^{\frac{t}{2}}\overline{F}(t-u)d G(u)}{\overline{F}(t)}
= m_{G}(\alpha)
- \alpha\int_{\frac{t}{2}}^{\infty}c(u)du
+ o\left( \int_{\frac{t}{2}}^{\infty}c(u)du +A_{2}(t) \right);
\end{eqnarray}

\item[(iii)] if $-2<\gamma<-1$,
\begin{equation}\label{eq3.26}
\frac{\int_{0}^{\frac{t}{2}}\overline{F}(t-u)d G(u)}{\overline{F}(t)}
= m_{G}(\alpha)
+ \alpha t c(t) \left( \int_{0}^{\frac{1}{2}}\left( (1-u)^{\beta}-1 \right)u^{\gamma}du +\frac{1}{2^{\gamma+1}(\gamma+1)} \right)
+ o\left( t c(t) +A_{2}(t) \right);
\end{equation}

\item[(iv)] if $\gamma=-2$ with $\int_{1}^{\infty} uc(u)du=\infty$,
\begin{eqnarray}\label{eq3.27}
\frac{\int_{0}^{\frac{t}{2}}\overline{F}(t-u)d G(u)}{\overline{F}(t)}
= m_{G}(\alpha)
- \alpha\beta t^{-1}\int_{1}^{\frac{t}{2}}u c(u)du
+ o\left( t^{-1}\int_{1}^{\frac{t}{2}}u c(u)du +A_{2}(t) \right);
\end{eqnarray}

\item[(v)] if $\gamma\leq -2$ with $\int_{1}^{\infty} uc(u)du<\infty$,
\begin{eqnarray}\label{eq3.28}
\frac{\int_{0}^{\frac{t}{2}}\overline{F}(t-u)d G(u)}{\overline{F}(t)}
= m_{G}(\alpha)
- \beta t^{-1}\int_{0}^{\infty}u e^{\alpha u}d G(u)
+ o\left( t^{-1} +A_{2}(t) \right).
\end{eqnarray}

\end{itemize}
\end{lemma}

\noindent {\bf Proof.}
(i). For $\gamma=-1$ with $\int_{0}^{\infty} e^{\alpha u}d G(u)=\infty$, we have
\begin{eqnarray*}
\lim_{t\to \infty} \frac{t c(t)}{\int_{0}^{t} c(u)du}=0
\end{eqnarray*}
by Lemma \ref{lem2.1}. Combining with Lemma \ref{lem2.3}, for large $t$ and arbitrary $\varepsilon>0$, we can get
\begin{eqnarray*}
& &
\frac{t c(t)}{\int_{0}^{\frac{t}{2}} c(y)dy}\left( \log 2 +\log u + \varepsilon\left( 2^{\varepsilon} -u^{-\varepsilon} \right) \right)
- \frac{c(t)}{\alpha\int_{0}^{\frac{t}{2}} c(y)dy} u^{-1}
\left( 1+\varepsilon^{2}u^{-\varepsilon} \right)
+\frac{1}{\alpha\int_{0}^{\frac{t}{2}} c(y)dy}\\
&\leq&
\frac{\C(t u)}{\alpha \int_{0}^{\frac{t}{2}} c(y)dy}-1 \\
&\leq&
\frac{t c(t)}{\int_{0}^{\frac{t}{2}} c(y)dy}\left( \log 2 +\log u - \varepsilon\left( 2^{\varepsilon} -u^{-\varepsilon} \right) \right)
- \frac{c(t)}{\alpha\int_{0}^{\frac{t}{2}} c(y)dy} u^{-1}
\left( 1-\varepsilon^{2}u^{-\varepsilon} \right)
+\frac{1}{\alpha\int_{0}^{\frac{t}{2}} c(y)dy},
\end{eqnarray*}
 if $1/t<u<1/2$. Then
\begin{eqnarray}\label{eq3.29}
& &
\beta \int_{\frac{1}{t}}^{\frac{1}{2}} (1-u)^{\beta-1}
\left( \frac{\C(tu)}{\alpha \int_{0}^{\frac{t}{2}} c(y)dy} -1 \right)du \nonumber \\
&=&
\frac{t c(t)}{\int_{0}^{\frac{t}{2}} c(y)dy}\left( \beta
\int_{0}^{\frac{1}{2}} (1-u)^{\beta-1}(\log u)du
+(\log 2)\left( 1-2^{-\beta} \right) + o(1) \right)
+ \frac{1-2^{-\beta}}{\alpha\int_{0}^{\frac{t}{2}} c(y)dy}(1+o(1))
\end{eqnarray}
for large $t$. Note that
\begin{eqnarray*}
\frac{\alpha\int_{0}^{\frac{t}{2}} c(y)dy}{\C\left( \frac{t}{2} \right)}
=1-\frac{1}{\alpha\int_{0}^{\frac{t}{2}}c(u)du}+ o\left( \frac{1+t c(t)}{\int_{0}^{\frac{t}{2}} c(y)dy} \right)
\end{eqnarray*}
for large $t$. Hence,
\begin{eqnarray}\label{eq3.30}
& &
\beta\int_{\frac{1}{t}}^{\frac{1}{2}} (1-u)^{\beta-1} \left( \frac{\C(tu)}{\C\left( \frac{t}{2} \right)} -1 \right) du \nonumber\\
&=&
\frac{\alpha \int_{0}^{\frac{t}{2}} c(y)dy}{\C\left( \frac{t}{2} \right)}
\beta \int_{\frac{1}{t}}^{\frac{1}{2}} (1-u)^{\beta-1}
\left( \frac{\C(tu)}{\alpha \int_{0}^{\frac{t}{2}} c(y)dy} -1 \right)du
+ \left( \frac{\alpha \int_{0}^{\frac{t}{2}} c(y)dy}{\C\left( \frac{t}{2} \right)}-1\right)
\beta \int_{\frac{1}{t}}^{\frac{1}{2}} (1-u)^{\beta-1}du \nonumber\\
&=&
\frac{t c(t)}{\int_{0}^{\frac{t}{2}} c(y)dy}\left( \beta
\int_{0}^{\frac{1}{2}} (1-u)^{\beta-1}(\log u)du
+(\log 2)\left( 1-2^{-\beta} \right)  \right)
+ o\left( \frac{1+t c(t)}{\int_{0}^{\frac{t}{2}} c(y)dy} \right).
\end{eqnarray}

Since
\begin{equation*}
\frac{\int_{1}^{\frac{t}{2}} \left( 1-\frac{u}{t} \right)^{\beta}
e^{\alpha u}d G(u)}{\C\left( \frac{t}{2} \right)}-1
=
\beta\int_{\frac{1}{t}}^{\frac{1}{2}} (1-u)^{\beta-1}
\left( \frac{\C(tu)}{\C\left( \frac{t}{2} \right)} -1 \right) du
-\frac{\beta}{t}(1+o(1))-\frac{\C(1)}{\C\left( \frac{t}{2} \right)}
+ \frac{\beta\C(1)}{t\C\left( \frac{t}{2} \right)}(1+o(1))
\end{equation*}
and
\begin{equation*}
\frac{\int_{0}^{1} \left( 1-\frac{u}{t} \right)^{\beta}
e^{\alpha u}d G(u)}{\C\left( \frac{t}{2} \right)}
= \frac{\C(1)}{\C\left( \frac{t}{2} \right)} -\frac{\beta}{t\C\left( \frac{t}{2}\right)}\int_{0}^{1}ue^{\alpha u}dG(u) (1+o(1))
\end{equation*}
for large $t$, we have
\begin{equation}\label{eq3.31}
\frac{\int_{0}^{\frac{t}{2}} \left( 1-\frac{u}{t} \right)^{\beta}
e^{\alpha u}d G(u)}{\C\left( \frac{t}{2} \right)}-1
=
\frac{t c(t)}{\int_{0}^{\frac{t}{2}} c(y)dy}\left( \beta
\int_{0}^{\frac{1}{2}} (1-u)^{\beta-1}(\log u)du
+(\log 2)\left( 1-2^{-\beta} \right)  \right)
+ o\left( \frac{1+t c(t)}{\int_{0}^{\frac{t}{2}} c(y)dy} \right)
\end{equation}
by combining with \eqref{eq3.30}.

Noting that
\[ \frac{\int_{0}^{\frac{t}{2}} \left( 1-\frac{u}{t} \right)^{\beta}
\frac{\left( 1-\frac{u}{t} \right)^{\rho_{2}}-1}{\rho_{2}}
e^{\alpha u}d G(u)}{\C\left( \frac{t}{2} \right)} \to 0\]
and
\[ \frac{\int_{0}^{\frac{t}{2}} \left( 1-\frac{u}{t} \right)^{\beta+\rho_{2}-\varepsilon}
e^{\alpha u}d G(u)}{\C\left( \frac{t}{2} \right)}\to 1 \]
as $t\to \infty$. From Lemma \ref{lem2.3} and \eqref{eq3.31}, it follows that
\begin{eqnarray*}
&&
\frac{\int_{0}^{\frac{t}{2}} \frac{b(t-u)}{b(u)}
e^{\alpha u}d G(u)}{\C\left( \frac{t}{2} \right)}\\
&=&
\frac{\int_{0}^{\frac{t}{2}} \left( 1-\frac{u}{t} \right)^{\beta}
e^{\alpha u}d G(u)}{\C\left( \frac{t}{2} \right)}
+ o(A_{2}(t)) \\
&=&
1+ \frac{t c(t)}{\int_{0}^{\frac{t}{2}} c(y)dy}\left( \beta
\int_{0}^{\frac{1}{2}} (1-u)^{\beta-1}(\log u)du
+(\log 2)\left( 1-2^{-\beta} \right)  \right)
+ o\left( \frac{1+t c(t)}{\int_{0}^{\frac{t}{2}} c(y)dy} +A_{2}(t) \right)
\end{eqnarray*}
for large $t$, which implies \eqref{eq3.24}.

(ii)-(v). Note that
\begin{equation*}
1-\beta x+ 2^{1-\beta}\beta(\beta-1)x^{2} \leq (1-x)^{\beta}
\leq 1-\beta x,\, 0\leq \beta \leq 1,
\end{equation*}
\begin{equation}\label{eq3.32}
1-\beta x \leq (1-x)^{\beta} \leq 1-\beta x+ 2^{1-\beta}\beta(\beta-1)x^{2}, \,
1\leq \beta \leq 2 \,\, \text{or}\,\, \beta <0,
\end{equation}
and
\begin{equation*}
1-\beta x \leq (1-x)^{\beta} \leq 1-\beta x+ 2^{-1}\beta(\beta-1)x^{2},\, \beta>2
\end{equation*}
for $0<x<\frac{1}{2}$. We only consider the case of $\beta <0$. For
the rest cases, the arguments are similarly and details are omitted
here.

By using integration by parts, we have
\begin{eqnarray}\label{eq3.33}
& & \int_{1}^{\frac{t}{2}} \left( 1-\frac{u}{t} \right)^{\beta}e^{\alpha u}d G(u)
-\int_{1}^{\infty} e^{\alpha u}d G(u) \nonumber\\
&=&
\left( \left( 1-t^{-1} \right)^{\beta} -1 \right)c(1)
- 2^{-\beta}c\left( \frac{t}{2} \right)
+ \alpha \int_{1}^{\frac{t}{2}} \left( \left( 1-\frac{u}{t} \right)^{\beta}-1  \right) c(u)du \nonumber\\
& &
- \alpha \int_{\frac{t}{2}}^{\infty} c(u)du
-\frac{\beta}{t}\int_{1}^{\frac{t}{2}}\left( 1-\frac{u}{t} \right)^{\beta-1}c(u)du
\end{eqnarray}
From Lemma \ref{lem2.1} and \eqref{eq3.32}, we can get
\begin{eqnarray*}
\int_{1}^{\frac{t}{2}}\left( \left( 1-\frac{u}{t} \right)^{\beta} -1 \right)c(u)du
= -\beta t^{-1}\int_{1}^{\infty}u c(u)du + o(t^{-1}),
\end{eqnarray*}
\begin{eqnarray}
\int_{1}^{\frac{t}{2}} \left( 1-\frac{u}{t} \right)^{\beta-1}c(u)du
= \int_{1}^{\infty}c(u)du+o(1)
\end{eqnarray}
and
\begin{eqnarray*}
\int_{\frac{t}{2}}^{\infty} c(u)du=o\left(t^{-1}\right)
\end{eqnarray*}
for large $t$, if $\gamma\leq -2$ with $\int_{1}^{\infty}u c(u)du<\infty$.
Then
\begin{eqnarray}\label{eq3.35}
& &
\int_{1}^{\frac{t}{2}} \left( 1-\frac{u}{t} \right)^{\beta}e^{\alpha u}d G(u)
-\int_{1}^{\infty} e^{\alpha u}d G(u) \nonumber\\
&=&
-\frac{\beta c(1)}{t} -\frac{\alpha \beta}{t} \int_{1}^{\infty} u c(u)du -\frac{\beta}{t}\int_{1}^{\infty} c(u)du + o\left(t^{-1}\right)\nonumber\\
&=&
-\frac{\beta}{t}\int_{1}^{\infty} u e^{\alpha u} d G(u) + o\left(t^{-1}\right)
\end{eqnarray}
for large $t$.

Since
\[ \int_{0}^{1} \left( \left( 1-\frac{u}{t} \right)^{\beta} -1 \right)e^{\alpha u}dG(u)
= -\frac{\beta}{t}\int_{0}^{1} ue^{\alpha u}dG(u)(1+o(1)),\]
\[0< \left| \int_{0}^{\frac{t}{2}} \left( 1-\frac{u}{t} \right)^{\beta}
\frac{\left( 1-\frac{u}{t} \right)^{\rho_{2}}-1}{\rho_{2}}e^{\alpha u}dG(u) \right|
<2^{-\beta}\frac{2^{-\rho_{2}}-1}{-\rho_{2}}m_{G}(\alpha)<\infty\]
and
\[ 0<\int_{0}^{\frac{t}{2}}
\left( 1-\frac{u}{t} \right)^{\beta+\rho_{2}-\varepsilon}e^{\alpha u}dG(u)
< 2^{-\beta-\rho_{2}+\varepsilon}m_{G}(\alpha)<\infty\]
for $\beta <0$, $\rho_{2} \leq 0$, we can get
\begin{eqnarray}\label{eq3.36}
\int_{0}^{\frac{t}{2}} \frac{\overline{F}(t-u)}{\overline{F}(t)}d G(u)
&=&
\int_{0}^{\frac{t}{2}} \frac{b(t-u)}{b(t)}e^{\alpha u}d G(u)\nonumber\\
&=&
\int_{0}^{\frac{t}{2}} \left( 1-\frac{u}{t} \right)^{\beta}e^{\alpha u}d G(u) + o(A_{2}(t)) \nonumber\\
&=&
m_{G}(\alpha)
+\int_{1}^{\frac{t}{2}} \left( 1-\frac{u}{t} \right)^{\beta}e^{\alpha u}d G(u)
-\int_{1}^{\infty}e^{\alpha u}d G(u) \nonumber\\
& &
-\frac{\beta}{t}\int_{0}^{1} u e^{\alpha u}d G(u) +o\left( t^{-1}+A_{2}(t) \right)\\
&=&
m_{G}(\alpha)-\beta t^{-1}\int_{0}^{\infty} u e^{\alpha u}d G(u) +o\left( t^{-1}+A_{2}(t) \right)\nonumber
\end{eqnarray}
by combining with \eqref{eq3.35}, which complete the proof of case (v).

For $\gamma=-2$ with $\int_{1}^{\infty} u c(u)du=\infty$, we have
\[\int_{1}^{\frac{t}{2}} \left( \left( 1-\frac{u}{t} \right)^{\beta} -1 \right) c(u)du
= -\beta t^{-1}\int_{1}^{\frac{t}{2}} u c(u)du(1+o(1))\]
and
\[\int_{\frac{t}{2}}^{\infty} c(u)du
= o\left( t^{-1}\int_{1}^{\frac{t}{2}}uc(u)du \right)\]
by using \eqref{eq3.32} and Lemma \ref{lem2.1}.
Combining \eqref{eq3.33} and \eqref{eq3.36}, we can get
\[
\int_{0}^{\frac{t}{2}} \frac{\overline{F}(t-u)}{\overline{F}(t)}d G(u)
=m_{G}(\alpha)- \alpha\beta t^{-1}\int_{1}^{\frac{t}{2}} u c(u)d u
+o\left( t^{-1}\int_{1}^{\frac{t}{2}} u c(u)d u+A_{2}(t) \right)
\]
for large $t$, which deduces the result in case (iv).

Similarly, we have
\begin{eqnarray*}
\int_{1}^{\frac{t}{2}} \left( 1-\frac{u}{t} \right)^{\beta}e^{\alpha u}d G(u)
-\int_{1}^{\infty} e^{\alpha u}d G(u)
=\alpha t c(t)\left( \int_{0}^{\frac{1}{2}} \left( (1-u)^{\beta} -1 \right)u^{\gamma}du +\frac{1}{2^{1+\gamma}(1+\gamma)} + o(1) \right)
\end{eqnarray*}
for $-2<\gamma<-1$, and
\begin{eqnarray*}
\int_{1}^{\frac{t}{2}} \left( 1-\frac{u}{t} \right)^{\beta}e^{\alpha u}d G(u)
-\int_{1}^{\infty} e^{\alpha u}d G(u)
=-\alpha \int_{\frac{t}{2}}^{\infty} c(u)du(1+o(1))
\end{eqnarray*}
for $\gamma =-1$ with $m_{G}(\alpha)<\infty$. By using \eqref{eq3.33} and
\eqref{eq3.36}, we can obtain \eqref{eq3.25} and \eqref{eq3.26}, respectively.

The proof is complete.
\qed

\noindent
{\bf Proof of Theorem \ref{th1}.}
For large $t$ and $A$ satisfied $0<A<f(t)$ and
$\chi'(A)<\alpha$, we have
\begin{eqnarray*}
 \frac{\int_{0}^{f(t)}\overline{F}(t-u)d F(u)}
{\int_{f(t)}^{\frac{t}{2}}\overline{F}(t-u)d F(u)}
&\leq&
\frac{e^{\chi(t)-2\chi\left(\frac{t}{2}\right)}\int_{0}^{f(t)} e^{\alpha u} d F(u)}
{\int_{f(t)}^{\frac{t}{2}} \frac{\overline{F}(t-u)}{\overline{F}^{2}\left( \frac{t}{2} \right)}d F(u)}\\
&\leq&
\frac{\left( \int_{0}^{A} e^{\alpha u}d F(u) + \left( \alpha -\chi'(A) \right) \left( f(t)-A \right)e^{\chi(f(t))}\right)e^{\chi(t)-2\chi\left(\frac{t}{2}\right)}}
{\int_{f(t)}^{\frac{t}{2}} \frac{\overline{F}(t-u)}{\overline{F}^{2}\left( \frac{t}{2} \right)}d F(u)}\\
&\leq&
\frac{\left( \int_{0}^{A} e^{\alpha u}d F(u) +  \alpha f(t)\right)e^{(1-\varepsilon)\left( 1-2^{1-\rho} \right)\chi(t)}}
{\int_{f(t)}^{\frac{t}{2}} \frac{\overline{F}(t-u)}{\overline{F}^{2}\left( \frac{t}{2} \right)}d F(u)}\\
&=&
o\left( A_{1}(t)\chi(t)+\frac{1}{\chi(t)}+ \frac{\chi^{\frac{1}{2}}(t)}{t} \right)
\end{eqnarray*}
by using Lemma \ref{lem3}, where $f(t)=\frac{1}{\chi'(t)}$. Combining Lemma \ref{lem3} and Lemma \ref{lem4},
we can get
\begin{eqnarray*}
& & \frac{\chi^{\frac{1}{2}}\left( \frac{t}{2} \right)\overline{F*F}(t) }
{t\overline{F}^{2}\left( \frac{t}{2} \right)}\\
&=&
\frac{2\chi^{\frac{1}{2}}\left( \frac{t}{2} \right)}
{t\overline{F}^{2}\left( \frac{t}{2} \right)}
\int_{0}^{\frac{t}{2}}\overline{F}(t-u)d F(u)
\left( 1+\frac{\overline{F}^{2}\left( \frac{t}{2} \right)}
{2\int_{0}^{\frac{t}{2}}\overline{F}(t-u)d F(u)} \right)\\
&=&
\frac{\chi^{\frac{1}{2}}\left( \frac{t}{2} \right)}
{\frac{t}{2}\overline{F}^{2}\left( \frac{t}{2} \right)}
\int_{f(t)}^{\frac{t}{2}}\overline{F}(t-u)d F(u)
\left( 1+\frac{\int_{0}^{f(t)}\overline{F}(t-u)d F(u)}
{\int_{f(t)}^{\frac{t}{2}}\overline{F}(t-u)d F(u)} \right)
\left( 1+\frac{\overline{F}^{2}\left( \frac{t}{2} \right)}
{2\int_{0}^{\frac{t}{2}}\overline{F}(t-u)d F(u)} \right)\\
&=&
\frac{\alpha}{2}\sqrt{\frac{\pi}{\rho(1-\rho)}}
-\frac{2^{\rho-5}\alpha\sqrt{\pi}(\rho-2)(\rho-3)}{(\rho(1-\rho))^{\frac{3}{2}}}
\frac{1}{\chi(t)}
-\frac{1}{2^{\rho}}\sqrt{\frac{\rho\pi}{1-\rho}}\frac{\chi(t)}{t}
+ o\left( A_{1}(t)\chi(t)+\frac{1}{\chi(t)}+ \frac{\chi^{\frac{1}{2}}(t)}{t} \right)
\end{eqnarray*}
for large $t$, which deduces the desired result.
\qed

\noindent
{\bf Proof of Proposition \ref{th3}.}
Combining Lemma \ref{lem7} and Lemma \ref{lem8}, we can derive the desired results.
\qed

\noindent
{\bf Proof of Theorem \ref{th4}.}
From \eqref{eq3.18}, it follows that
\[ \overline{F}\left( \frac{t}{2} \right)\overline{G}\left( \frac{t}{2} \right)
=2^{-\beta -\gamma}\left( 1+A_{2}(t)\frac{2^{-\rho_{2}}-1}{\rho_{2}}
+A_{3}(t)\frac{2^{-\rho_{3}}-1}{\rho_{3}} + o\left(A_{2}(t)+A_{3}(t)
\right)\right)\overline{G}(t)b(t)\] for large $t$. By arguments
similar to Proposition \ref{th3}, we have
\[\int_{0}^{\frac{t}{2}}\overline{G}(t-u)dF(u)
+ \overline{F}\left( \frac{t}{2} \right)
\overline{G}\left( \frac{t}{2} \right)
=\overline{G}(t)M_{2}(\beta,\gamma,t)\]
for large $t$, where $M_{2}(\beta,\gamma,t)$ is given by \eqref{eq2.17} and \eqref{eq2.18}.
With the decomposition of convolution tail, for large $t$ we can get
\begin{eqnarray*}
\overline{F*G}(t)
&=& \int_{0}^{\frac{t}{2}}\overline{F}(t-u)d G(u)  + \int_{0}^{\frac{t}{2}}\overline{G}(t-u)d F(u)
+ \overline{F}\left( \frac{t}{2} \right)\overline{G}\left( \frac{t}{2} \right) \\
&=&
\overline{F}(t)M_{1}(\beta,\gamma,t)+ \overline{G}(t)M_{2}(\beta,\gamma,t)
\end{eqnarray*}
with $M_{1}(\beta,\gamma,t)$ given by \eqref{eq2.7} and \eqref{eq2.11}.
The proof is complete.
\qed

\vspace{1cm}

\noindent {\bf Acknowledgements}~~This work was supported by the National Natural Science Foundation
of China grant no.11171275, the Natural Science Foundation Project of CQ no. cstc2012jjA00029.

\end{document}